\documentclass[a4paper]{article}

\usepackage[english]{babel}
\usepackage[utf8]{inputenc}
\usepackage{amssymb,amsmath}
\usepackage{graphicx}
\usepackage[colorinlistoftodos]{todonotes}
\usepackage{algorithm,algorithmic}
\usepackage{multicol,adjustbox}
\usepackage[font=small,labelfont=bf,tableposition=top]{caption}
\usepackage{url}

\DeclareCaptionLabelFormat{andtable}{#1~#2  \&  \tablename~\thetable}

\DeclareMathOperator{\psd}{psd}
\DeclareMathOperator{\rank}{rank}
\DeclareMathOperator{\rankp}{\rank_{\psd}}
\DeclareMathOperator{\xc}{xc}
\DeclareMathOperator{\xcp}{\xc_{\psd}}
\DeclareMathOperator{\trace}{tr}

\newtheorem{conj}{Conjecture}
\newtheorem{observation}{Observation}
\newtheorem{exmp}{Example}[section]

\title{Algorithms for Positive Semidefinite Factorization}

\author{%
  Arnaud Vandaele\thanks{Department of Mathematics and Operational Research, Facult\'e Polytechnique, Universit\'e de Mons, Rue de Houdain~9, 7000 Mons, Belgium. Email: \texttt{\{arnaud.vandaele,nicolas.gillis\}@umons.ac.be}. NG acknowledges the support by the
F.R.S.-FNRS (incentive grant for scientific research no F.4501.16) and by the ERC (starting grant no 679515).}
  \and
  François Glineur\thanks{Center for Operations Research and Econometrics, Universit\'e catholique de Louvain, Voie du Roman Pays, 34, B-1348 Louvain-La-Neuve, Belgium;   
              ICTEAM Institute, Universit\'e catholique de Louvain, B-1348 Louvain-La-Neuve, Belgium.
              Email: \texttt{francois.glineur@uclouvain.be}. This paper presents research results of the Concerted Research Action (ARC) programme supported by the Federation Wallonia-Brussels (contract ARC 14/19-060). }
  \and
  Nicolas Gillis\footnotemark[1]
}

\date{}

\begin{document}
\maketitle

\begin{abstract}
This paper considers the problem of positive semidefinite factorization (PSD factorization), a generalization of exact nonnegative matrix factorization. 
Given an $m$-by-$n$ nonnegative matrix $X$ and an integer $k$, the PSD factorization problem consists in finding, 
if possible, symmetric $k$-by-$k$ positive semidefinite matrices $\{A^1,...,A^m\}$ and $\{B^1,...,B^n\}$ such that 
$X_{i,j}=\text{trace}(A^iB^j)$ for $i=1,...,m$, and $j=1,...,n$. 
PSD factorization is NP-hard. In this work, we introduce several local optimization schemes to tackle this problem: a fast projected gradient method and two algorithms based on the coordinate descent framework.  
The main application of PSD factorization is the computation of semidefinite extensions, that is, the representations of polyhedrons as projections of spectrahedra, for which the matrix to be factorized is the slack matrix of the polyhedron. We compare the performance of our algorithms on this class of problems. 
In particular, we compute the PSD extensions of size $k=1+ \lceil \log_2(n) \rceil$ for the regular $n$-gons when $n=5$, $8$ and $10$.  
We also show how to generalize our algorithms to compute the square root rank (which is the size of the factors in a PSD factorization where all factor matrices $A^i$ and $B^j$ have rank one) and completely PSD factorizations (which is the special case where the input matrix is symmetric and equality $A^i=B^i$ is required for all~$i$).   
\end{abstract}


\section{Introduction} \label{sec-intro} 

Given an $m$-by-$n$ nonnegative matrix $X$ and an integer $k<\min\{m,n\}$, the standard nonnegative matrix factorization (NMF) problem seeks a  matrix $\tilde{X}$ such that the $(i,j)$th entry of $\tilde{X}$ is equal the inner product of two size-$k$ nonnegative vectors $w_i$ and $h_j$, and which is as close to $X$ as possible. For all $i=1,...,m$ and $j=1,...,n$, we have 
\begin{equation} \label{nmf-entrywise}
	X_{ij}\approx\tilde{X}_{ij}=\langle w_i,h_j\rangle=w_i^Th_j \text{ with } w_i \text{ and } h_j \in \mathbb{R}^k_+.
\end{equation} 
Considering the set $\{w_1,...,w_m\}$ as the rows of a matrix $W$ and the set $\{h_1,...,h_n\}$ as the columns of a matrix $H$, the approximating matrix $\tilde{X}$ is the product of two nonnegative matrices $W\in\mathbb{R}^{m\times k}$ and $H\in\mathbb{R}^{k\times n}$.
This leads to the following optimization problem for NMF: 
\begin{equation} \label{nmf-opt}
	\min_{W\geq 0,H\geq 0}\|X-WH\|_F^2, 
\end{equation}
where $||X||_F^2 = \sum_{i,j} X_{ij}^2$ is the Frobenius norm of matrix $X$. 
NMF (\ref{nmf-opt}) has become a widely-used approach for linear dimensionality reduction. 
In fact, when the columns of the matrix $X$ represent the elements of a data set, the nonnegative factorization allows to interpret the $j$th column of $X$ as a nonnegative linear combination of the columns of $W$ where the weights are given by the $j$th column of $H$. 
Unlike other comparable techniques, the nonnegativity imposed on the entries of $W$ and $H$ leads to a better interpretability of the decomposition and NMF 
has been proved successful in many fields of data analysis such as image processing, text mining and hyperspectral imaging; see~\cite{G14} and the references therein. 

The problem of PSD factorization addressed in this paper is a recently introduced generalization of NMF~\cite{GPT13}. In a PSD factorization problem, the cone of positive semidefinite matrices replaces the nonnegative orthant of NMF.
More precisely, the inputs of a PSD factorization problem are the same as for NMF, namely, 
a $m$-by-$n$ nonnegative matrix $X$ and an integer $k<\min\{m,n\}$.
However, instead of using the inner product between two vectors of size $k$, the $(i,j)$th entry of the approximating matrix $\tilde{X}$ is given by the inner product between two symmetric $k$-by-$k$ positive semidefinite matrices, $A^i$ and $B^j$.
The inner product of two matrices is a generalization of the dot product of two vectors, and is equal to the trace of the product of the two matrices. 
Hence, we have for $i=1,...,m$ and $j=1,...,n$, 
\begin{equation*} 
	X_{ij}\approx\tilde{X}_{ij}=\langle A^i,B^j\rangle=\trace(A^iB^j) \text{ with } A^i \text{ and } B^j \in \mathcal{S}_+^{k}.
\end{equation*}

As in NMF, the optimization problem corresponding to PSD factorization consists in minimizing the quantity $\|X-\tilde{X}\|^2_F$.
It can be expressed by the following non-convex and NP-hard problem~\cite{shitov2016complexity} where the variables are the two sets of matrices $\{A^1,...,A^m\}$ and $\{B^1,...,B^n\}$ belonging to the positive semidefinite cone $\mathcal{S}_+^{k}$: 
\begin{equation} \label{psd-opt}
	\min_{\substack{A^i,B^j\in\mathcal{S}_+^{k} \\ i=1,...,m \\ j=1,...,n}}\sum_{i=1}^m\sum_{j=1}^n\left(X_{ij}-\left\langle A^i,B^j\right\rangle\right)^2.
\end{equation}

In this work, we propose several algorithms for solving (\ref{psd-opt}) numerically.
To our knowledge, no algorithm has been proposed in the literature to solve this problem. 
The paper is organized as follows. In Section~\ref{sec-defmot}, we introduce the PSD factorization problem in more details, highlighting its connection with extended formulations.  
In Section~\ref{sec-results}, we propose several algorithms to compute PSD factorizations (namely, a fast projected gradient
method and two algorithms based on the coordinate descent framework). 
In Section~\ref{sec-results}, we compare the efficiency of the presented methods on a benchmark of nonnegative matrices. 
In Section~\ref{sec-application}, we show how to use our algorithms to compute 
(i)~PSD factorizations of the slack matrices of regular $n$-gons, 
(ii)~symmetric PSD factorizations related to completely PSD matrices, and 
(iii)~the square root rank of nonnegative matrices.

\section{Linear and semidefinite extensions, and factorizations} \label{sec-defmot}



\begin{figure}[h!]
	\begin{center}
		\includegraphics[scale=0.2]{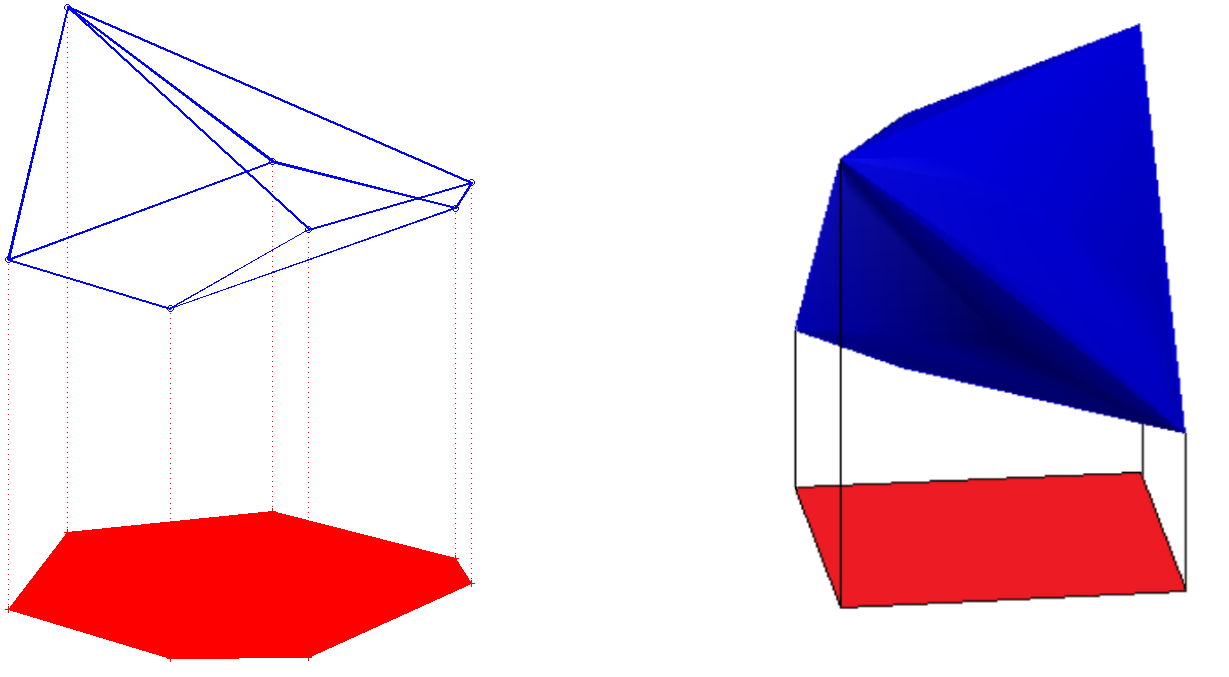}
		\caption{Left: linear extension of size $6$ of an irregular hexagon. Right: psd-lift of size $3$ of the square.}
		\label{extensions}
	\end{center}
\end{figure}

In the context of the NMF of an $m$-by-$n$ nonnegative matrix $X$, the minimum value of the inner dimension $k$ for which it is possible to find $W\in\mathbb{R}^{m\times k}_+$ and $H\in\mathbb{R}^{k\times n}_+$ such that $X=WH$ is called the nonnegative rank of $X$, and is denoted $\rank_+(X)$.
The search for such an exact factorization has tight connections with the study of linear extensions of polyhedrons.
Let $P$ be a polyhedron described by a system of linear inequalities.
A linear extension of $P$ is another polyhedron $Q$ of higher dimension which projects linearly onto $P$, that is, for which there exists a linear map $\pi$ such that $\pi(Q) = P$.
Such linear extensions are particularly useful when the size (measured by the number of facets) of a linear extension is (much) smaller than the size of the initial polyhedron.
For example, the left picture of Figure \ref{extensions} illustrates a linear extension of an irregular (planar) heptagon, which is three-dimensional but features only six facets.
Among all possible linear extensions of a polyhedron $P$, the size of the smallest one is the \textit{linear extension complexity} of $P$ and is denoted by $\xc(P)$.

An outstanding result of Yannakakis establishes a strong link between NMF and linear extensions \cite{Y91}: the linear extension complexity of a polyhedron $P$ is equal to the nonnegative rank of a particular matrix related to $P$, called the \emph{slack matrix} $\mathcal{S}_P$: 
\begin{equation}
	\xc(P)=\rank_+(\mathcal{S}_{P}). \label{yannakakis}
\end{equation}
For a polyhedron $P$ featuring $f$ facets and $v$ vertices, the slack matrix $\mathcal{S}_{P}$ is a $f$-by-$v$ nonnegative matrix whose $(i,j)$th entry is the slack between the $i$th facet and the $j$th vertex. Furthermore, Yannakakis showeed that any rank-$k$ nonnegative factorization of $\mathcal{S}_{P}$ (implicitly) provides a size-$k$ linear extension of $P$. 
This result connecting the two fields has been at the core of many recent developments; see, e.g.,~\cite{K11} and the references therein. 

Recently, the work of Yannakakis was generalized to allow for arbitrary closed convex cones $K$ instead of the nonnegative orthant~\cite{GPT13}. From the point of view of extensions, linear extensions (for which $K=\mathbb{R}_+^k$) are replaced by conic extensions, that is, representations as projections under a linear map of an affine slice of a cone $K$. These generalized extensions are also called \textit{$K$-lifts}.
In this paper, we focus on the case where $K=\mathcal{S}^k_+$, the cone of positive semidefinite matrices. 
In that case, given a polyhedron $P$, we are looking for a \emph{spectrahedron} (an affine slice of a positive semidefinite cone) which projects onto $P$ under a linear map. Moreover, we are trying to find such a \emph{semidefinite extension} whose size (as measured by the dimension of the positive semidefinite cone) is as small as possible. This minimal size is called the \emph{semidefinite extension complexity} of $P$, and is denoted $\xcp(P)$ .

The semidefinite extension complexity never exceeds the linear extension complexity, but can be strictly lower. For example, the linear extension complexity of the square is $4$, but there exists a spectrahedron of size 3 which projects linearly onto the square (see the picture on the right of Figure \ref{extensions}).
Yannakakis' result (\ref{yannakakis}) can be generalized in the following way, which uses the \textit{positive semidefinite rank} (abbreviated psd-rank or $\rankp$) 
to a special rank of the slack matrix of $P$~\cite[Theorem 3.3]{GPT13}:  
\begin{equation*} 
	\xcp(P)=\rankp(\mathcal{S}_P).
\end{equation*}
The {positive semidefinite rank} is related to the PSD factorization problem (\ref{psd-opt}) in the same way than the nonnegative rank is connected to NMF. 
Formally, the psd-rank of a $m$-by-$n$ nonnegative matrix $X$ is the smallest integer $k$ for which there exist two sets of $k$-by-$k$ positive semidefinite matrices $\{A^1,...,A^m\}$ and $\{B^1,...,B^n\}$ such that $X_{ij}=\langle A^i,B^j \rangle$ holds for all $i=1,...,m$ and $j=1,...,n$.
We refer the reader to the survey \cite{FGP14} for further informations on the psd-rank.

\begin{exmp} \label{exemplepsdfactorization}
In order to illustrate the concept of the size of a PSD factorization, let the following $4$-by-$4$ matrix be a slack matrix of the square, 
$$S_4=\begin{pmatrix}
0 & 1 & 1 & 0 \\ 0 & 0 & 1 & 1 \\ 1 & 0 & 0 & 1 \\ 1 & 1 & 0 & 0
\end{pmatrix}.$$
Already highlighted by the picture on the right of Figure \ref{extensions}, it is possible to find a $\mathcal{S}_+^3$ factorization of $S_4$, for example with the following factors:
$$
 A^1=\begin{pmatrix} 1 & 0 & 0 \\  0 & 0 & 0 \\  0 & 0 & 0 \end{pmatrix}, A^2=\begin{pmatrix} 0 & 0 & 0 \\  0 & 1 & 0 \\  0 & 0 & 0 \end{pmatrix}, A^3=\begin{pmatrix} 0 & 0 & 0 \\  0 & 0 & 0 \\  0 & 0 & 1 \end{pmatrix}, A^4=\begin{pmatrix} 1 & -1 & 1 \\  -1 & 1 & -1 \\  1 & -1 & 1 \end{pmatrix},$$
$$B^1=\begin{pmatrix} 0 & 0 & 0 \\  0 & 0 & 0 \\  0 & 0 & 1 \end{pmatrix}, B^2=\begin{pmatrix} 1 & 0 & 0 \\  0 & 0 & 0 \\  0 & 0 & 0 \end{pmatrix}, B^3=\begin{pmatrix} 1 & 1 & 0 \\  1 & 1 & 0 \\  0 & 0 & 0 \end{pmatrix}, B^4=\begin{pmatrix} 0 & 0 & 0 \\  0 & 1 & 1 \\  0 & 1 & 1 \end{pmatrix}.
$$
\end{exmp}
Designing algorithms for solving (\ref{psd-opt}) is therefore of great interest in the search of psd-lifts based on the factorization of the corresponding slack matrices, and is the main objective of this paper.

\section{Algorithms for PSD factorization} \label{sec-algos}

The PSD factorization problem (\ref{psd-opt}) is nonconvex.
However, when one of the two sets of matrix variables $\{A^1,...,A^m\}$ or $\{B^1,...,B^n\}$ is fixed, optimizing over the other set reduces to a convex problem.
For this reason, we develop in this work algorithms using an alternating strategy for solving (\ref{psd-opt}), by optimizing alternately over the sets $\{A^1,...,A^m\}$ and $\{B^1,...,B^n\}$.
The same approach is used by nearly all NMF algorithms for solving (\ref{nmf-opt}), which is also a nonconvex problem 
that becomes convex when one of the two factors is fixed. 
The pseudo-code of the general alternating scheme for PSD factorization is detailed in Algorithm~\ref{algpsd-alternate}. 


\begin{algorithm}[h!]
\caption{Alternating Strategy for PSD Factorization}
\label{algpsd-alternate}
\begin{algorithmic}[1]
\STATE INPUT: $X\in\mathbb{R}_+^{m\times n}$ and initial iterates $\{A^1,...,A^m\}$ and $\{B^1,...,B^n\}$.
\STATE OUTPUT: $\{A^1,...,A^m\}$ and $\{B^1,...,B^n\}$.
\WHILE{stopping criterion not satisfied}
\STATE $\{A^1,...,A^m\}\leftarrow \text{optimize subproblem}(X,\{B^1,...,B^n\})$, \label{psd-alter-3}
\STATE $\{B^1,...,B^n\}\leftarrow \text{optimize subproblem}(X^T,\{A^1,...,A^m\})$. \label{psd-alter-4}
\ENDWHILE
\end{algorithmic}
\end{algorithm}

Since the subproblems are symmetric, we can assume without loss of generality for the presentation of the algorithms that the set $\{B^1,...,B^n\}$ is fixed and that we want to optimize over the $A^i$'s. 
The corresponding problem can be written formally as:

\begin{equation}\label{psd-obj-sub}
\min_{\substack{A^i\in\mathcal{S}_+^{k} \\ i=1,...,m}}\sum_{i=1}^m\sum_{j=1}^n\left(X_{ij}-\left\langle A^i,B^j\right\rangle\right)^2.
\end{equation}

Matrices $A^i$ do not influence each other in (\ref{psd-obj-sub}), that is, the problem is separable, hence it reduces to $m$ independent convex problems, each corresponding to the optimization over a single factor $A^i$ (corresponding to a single row of $X$). 
Hence, our first idea consists in solving each of these problems to optimality as described by Algorithm~\ref{algsub-exact}, which is an instance of a semidefinite program.
The combination of Algorithms \ref{algpsd-alternate} and \ref{algsub-exact} leads to an exact two-block coordinate descent scheme.
Since each block of variables belong to a closed convex set and the objective function is continuously differentiable, a stationary point of (\ref{psd-opt}) is obtained in the limit~\cite{grippo2000convergence}. 

\begin{algorithm}[h!]
\caption{optimize subproblem (exact)}
\label{algsub-exact}
\begin{algorithmic}[1]
\STATE INPUT: $X\in\mathbb{R}_+^{m\times n}$ and $B^j\in\mathcal{S}_+^{k}$ with $j=1,...,n$.
\STATE OUTPUT: $\{A^1,...,A^m\}$
\FOR{$i=1$ \TO $m$}
\STATE $A^i\leftarrow \arg \min_{\substack{A^i\in\mathcal{S}_+^{k}}} \sum_{j=1}^n\left(X_{ij}-\left\langle A^i,B^j\right\rangle\right)^2$
\ENDFOR
\end{algorithmic}
\end{algorithm}

We implemented Algorithm~\ref{algsub-exact} with the general convex solver YALMIP \cite{lofberg2004yalmip}.
However, this approach has proven to be far too slow in comparison with the other methods developed hereafter.
As is also the case in the context of NMF, 
the reason for the poor performance is that it is not worth solving the subproblems (\ref{psd-obj-sub}) to optimality at each iteration. Once the objective function has decreased by some amount, it is preferable to move quickly to the other set of variables rather than performing extra work to refine the subproblem solution to optimality. 
Based on that observation, we propose in the following two iterative methods for solving (\ref{psd-obj-sub}): an algorithm based on the (accelerated) gradient method described in Section~\ref{sec-gradient}, 
and implementations of coordinate descent methods introduced in Section~\ref{sec-cd}.

\subsection{A fast projected gradient method} \label{sec-gradient}


One of the most widely used method in continuous optimization is the \textit{gradient method}.
From a given starting point $x_0$, a sequence of points $\{x_t\}$ is built by taking a step in the direction $-\nabla f(x_{t-1})$ for each iterate $t=1,2,...$. 
The next point is then computed as $x_{t}=x_{t-1}-\alpha_{t-1}\nabla f(x_{t-1})$, where quantity $\alpha_{t-1}$ is the step size along the steepest descent direction.
The gradient method admits accelerated schemes, which were first introduced in~\cite{nesterov1983method}.
The scheme used in this work is described as Algorithm~\ref{algfgradientstep} for the general problem $\min_{x\in Q} f(x)$ with $Q$ a closed convex set.

\begin{algorithm}[h!]
\caption{Nesterov's accelerated gradient method}
\label{algfgradientstep}
\begin{algorithmic}[1]
\STATE INPUT: $x_0\in Q$
\STATE OUTPUT: $x_{\text{maxiter}}$
\STATE Set $x_{-1}=x_0$
\FOR{$t=1:\text{maxiter}$}
\STATE $y_{t}=x_{t-1}+\frac{t-2}{t+1}(x_{t-1}-x_{t-2})$ \label{alg-fgradientstep-line5}
\STATE $x_{t}=\text{Proj}_Q(y_{t}-\frac{1}{L}\nabla f(y_{t}))$ \label{alg-fgradientstep-line6}
\ENDFOR
\end{algorithmic}
\end{algorithm}

The accelerated gradient method presented as Algorithm~\ref{algfgradientstep} has roughly the same computational cost as the usual gradient method. 
The difference lies in the fact that the gradient step (line \ref{alg-fgradientstep-line6}) is made at an extrapolation point $y_t$ (computed in line \ref{alg-fgradientstep-line5}) instead of the previous iterate $x_{t-1}$.
When using a step size equal to $\frac{1}{L}$ (with $L$ the Lipschitz constant of the objective function's gradient $\nabla f$), the accelerated gradient method exhibits a convergence rate of $O(1/t^2)$, with $t$ the number of iterations (see \cite{nesterov1983method} for more details).
In order to apply the accelerated scheme of Algorithm~\ref{algfgradientstep} to the PSD factorization problem (\ref{psd-obj-sub}), several issues must first be addressed.

\begin{itemize}
	\item \textbf{Computing the gradient.} Let $f$ denote the quantity to minimize in~(\ref{psd-obj-sub}).
Using the Frobenius norm, $f$ can be written as follows, 
	\begin{equation} \label{psd-fctaderiver}
		f = \sum_{i=1}^m \sum_{j=1}^n \left(X_{ij}-\left\langle A^i,B^j \right\rangle\right)^2  
		= \left\|X - \mathcal{A}^T\mathcal{B}\right\|^2_F, 
	\end{equation}
with $\mathcal{A}=\begin{pmatrix}vec(A^1) \hspace{0.1cm} ... \hspace{0.1cm} vec(A^m)\end{pmatrix}$ and $\mathcal{B}=\begin{pmatrix}vec(B^1) \hspace{0.1cm} ... \hspace{0.1cm} vec(B^n)\end{pmatrix}$ being $k^2$-by-$m$ and $k^2$-by-$n$ matrices respectively.
Using this notation, the gradient of $f$ with respect to the variable $\mathcal{A}$ is:
	\begin{equation*}
		\nabla f = -2\left(X-\mathcal{A}^T\mathcal{B}\right)\mathcal{B}^T.
	\end{equation*}
From (\ref{psd-fctaderiver}), we can also derive the Lipschitz constant $L$ needed in Algorithm~\ref{algfgradientstep}, which will be equal to the largest eigenvalue of the Hessian $\nabla^2 f$, hence  $L=2\lambda_{\max}\left(\mathcal{B}\mathcal{B}^T\right)$. 

    \item \textbf{Projecting on $\mathcal{S}^k_+$.}
    For our problem, the closed convex set $Q$ that we need to project onto (see line \ref{alg-fgradientstep-line6} of Algorithm~\ref{algfgradientstep}) is the cone of symmetric and positive semidefinite matrices, that is, $Q=\mathcal{S}_+^{k}$.
    For every $k$-by-$k$ real symmetric matrix $C$, we have $C=U\Lambda U^{T}$ where $U$ is an orthogonal matrix, and $\Lambda$ is a diagonal matrix whose entries are the eigenvalues of $C$.
    Defining $\Lambda_+=\text{diag}\left(\max(0,\lambda_1),...,\max(0,\lambda_k)\right)$, we have the following closed-form formula for the projection
	\begin{equation} \label{psd-proj-opti}
	\text{Proj}_{\mathcal{S}^k_+}\left(C\right) = \arg\min_{X\in\mathcal{S}^k_+}\|X-C\| = U\Lambda_+ U^{T}.
	\end{equation}
The main computational cost of the projection (\ref{psd-proj-opti}) is spent computing the spectral decomposition of~$C$.
\end{itemize}

The pseudo-code of the accelerated gradient method for PSD factorization is presented as Algorithm~\ref{algfgradient} and denoted FPGM (for Fast Projected Gradient Method).
Recall that Algorithm~\ref{algfgradient} is used for solving the subproblems of the general alternating scheme of Algorithm~\ref{algpsd-alternate}.
We choose to perform a (fixed) number of accelerated gradient steps proportional to the size of the factors, equal to $k \Delta$ where $\Delta$ is a parameter   
(line \ref{nbsteps} of Algorithm~\ref{algfgradient}).
In Section \ref{sec-results}, performance of the algorithm is compared for different values of $\Delta$.

\begin{algorithm}[h!]
\caption{optimizesubproblem (FPGM)}
\label{algfgradient}
\begin{algorithmic}[1]
\STATE INPUT: $X\in\mathbb{R}_+^{m\times n}$ and $B^j\in\mathcal{S}_+^{k}$ with $j=1,...,n$, parameter $\Delta$. 
\STATE OUTPUT: $\{A^1,...,A^m\}$
\STATE $\{A^1,...,A^m\} \leftarrow$ Initialization
\STATE Construct $\mathcal{A}_0$ from $\{A^1,...,A^m\}$ and $\mathcal{B}$ from $\{B^1,...,B^n\}$.
\STATE $L\leftarrow \lambda_{\max}\left(\mathcal{B}\mathcal{B}^T\right)$
\STATE Set $\mathcal{A}_{-1}=\mathcal{A}_0$
\FOR{$t=1:k\Delta$} \label{nbsteps}
\STATE $\mathcal{Y}_{t}=\mathcal{A}_{t-1}+\frac{t-2}{t+1}(\mathcal{A}_{t-1}-\mathcal{A}_{t-2})$
\STATE $\mathcal{A}_{t}=\text{Proj}_{\mathcal{S}_+^k}\left(\mathcal{Y}_t+\frac{1}{L}(X\mathcal{B}^T-\mathcal{Y}_t^T\mathcal{B}\mathcal{B}^T)\right)$ \label{alg-fgradient-line6}
\ENDFOR
\STATE Extract $\{A^1,...,A^m\}$ from $\mathcal{A}_{k\Delta}$. 
\end{algorithmic}
\end{algorithm}


Algorithm~\ref{algfgradient} has two drawbacks. First, in the current form of the algorithm, it is not possible to adjust easily the rank of the $A^i$'s and the $B^j$'s while it is interesting to obtain low-rank factors (observe that the factors of Example~\ref{exemplepsdfactorization} are all rank one); 
see the discussion in Section~\ref{innerrank}.  
Second, if we know beforehand the values of some entries of the $A^i$'s and the $B^j$'s, it is not straightforward to keep them constant during the iterations of the algorithm (the projection step would become even more computationally expensive, as a linearly constrained semidefinite program would have to be solved). 
In the next section, we present coordinate-descent algorithms overcoming these limitations.

\subsection{Coordinate descent algorithms} \label{sec-cd}
Although known for many yeras, coordinate descent (CD) methods have recently received a new lease of life \cite{wright2015coordinate}.
This increase of interest is mainly due to the increasing number of large-scale optimization problems in data mining and machine learning applications for which the simplicity of the CD methodology allows efficient and competitive implementations (while high solution accuracy is usually not needed since data is typically rather noisy).   
In many of these applications, fixing all variables except one leads to an univariate optimization problem for which computation of a minimizer is cheap. 
For example, in the case of the NMF problem, the corresponding univariate optimization problem is quadratic, and its optimal solution can therefore be written in closed form.
First introduced in \cite{cichocki2007hierarchical} under the name HALS (for Hierarchical Alternating Least Square), methods based on the CD scheme have proven to be among the most effective ones for the NMF problem~\cite{CP09b, hsieh2011fast, GG12}. 

\subsubsection{Change of variables}

If we want to successfully apply the CD scheme to the PSD factorization problem, it is crucial that the update of one variable is  computationally cheap and easy to implement. 
However, it is not straightforward to update the entries of the factors $A^i$ and $B^j$: unlike NMF where nonnegativity of the variables had to be taken into account, which can be ensured separately in each variable (a separable constraint), matrices $A^i$ and $B^j$ are required to remain positive semidefinite, which is no longer separable.
Hence, in order to adapt the problem (\ref{psd-obj-sub}) to the application of the CD scheme, we perform a simple change of variables popularized by the works of Burer and Monteiro on semidefinite programming~\cite{burer2003nonlinear}. 
Since every symmetric positive semidefinite matrix $L$ can be written in the form $L=HH^T$, we introduce new (matrix) variables $a^i\in\mathbb{R}^{k\times r_i}$ for $i=1,...,m$ and $b^j\in\mathbb{R}^{k\times r_j}$ for $j=1,...,n$ linked to the original factors $A^i$ and $B^j$ as follows: 
\begin{equation*} 
	A^i=a^i{a^i}^T \text{ and } B^j=b^j{b^j}^T.
\end{equation*}  
With this reformulation, entries of the new variables $a^i$ and $b^j$ are unconstrained, and poitive semidefiniteness of the $A^i$'s and $B^j$'s  is automatically guaranteed.  
Another benefit is the ability to easily adjust the inner rank of the factors $A^i$ and $B^j$ by choosing the number of columns $r$ of the new variables, as the rank will be at most equal to this number.
Moreover, if some entries of $a^i$ or $b^j$ are known and fixed, they can simply be ignored in the CD scheme.
Since we have $\left\langle A^i,B^j\right\rangle=\sum_{h=1}^{r_i}\sum_{l=1}^{r_j}\left({a_{:,h}^i}^Tb_{:,l}^j\right)^2$, the optimization problem (\ref{psd-obj-sub}) is now written as follows with the new variables:
\begin{equation} \label{psd-obj-newvariables}
\min_{a^i\in\mathbb{R}^{k\times r_i}, i=1,...,m } f=\sum_{i=1}^m\sum_{j=1}^n\left(X_{i,j}-\sum_{h=1}^{r_i}\sum_{l=1}^ {r_j}\left({a_{:,h}^i}^Tb_{:,l}^j\right)^2\right)^2.
\end{equation}
Since we use an alternating scheme, we assume in the remainder of the section that the $b^j$'s are given and that only the $a^i$'s must be optimized. Note that the matrix $a^i \in\mathbb{R}^{k\times r_i}$ is made of $kr_i$ entries, so that the number of variables of the problem (\ref{psd-obj-newvariables}) is $k\sum_{i=1}^mr_i$, and $mk^2$ in the full-rank case ($r_i = k$ for all $i$). 

\subsubsection{Update of one variable}
In order to apply the CD scheme to (\ref{psd-obj-newvariables}), we need to derive the expression of the univariate function to minimize when all the variables of (\ref{psd-obj-newvariables}) are fixed but one, say the entry $(p,q)$ of the factor $a^i$ denoted $a_{p,q}^i$. 
Only the $i$th factor is impacted when the entry $a_{p,q}^i$ is updated since the factors $A^i$'s are independent from one another. 
By highlighting $a_{p,q}^i$, the part of the objective function influenced by the variable is
\footnotesize
\begin{equation} \label{psd-contribution}
\sum_{j=1}^n\left(X_{i,j}-\sum_{\substack{h=1 \\ h\neq q}}^r\sum_{l=1}^r\left({a_{:,h}^i}^Tb_{:,l}^j\right)^2-\sum_{l=1}^r\left({a_{\bar{p},q}^i}^Tb_{\bar{p},l}^j\right)^2-\textcolor{red}{ \boxed{a_{p,q}^i} }\left(2{a_{\bar{p},q}^i}^T\left(\sum_{l=1}^rb_{p,l}^jb_{\bar{p},l}^j\right)\right)-\textcolor{red}{ \boxed{{a_{p,q}^i}^2} }\|b_{p,:}^j\|^2\right)^2,
\end{equation}
\normalsize
where $\bar{p}=\{1,...,k\}\backslash \{p\}$.
We observe that the function to minimize is a fourth degree polynomial in $a_{p,q}^i$ and its gradient has therefore the form of a cubic polynomial,
\begin{equation} \label{psd-gradient}
\nabla_{a_{p,q}^i}f = c_3{a_{p,q}^i}^3+c_2{a_{p,q}^i}^2+c_1{a_{p,q}^i}+c_0,
\end{equation}
where
\small
\begin{eqnarray*}
c_3 \hspace{-0.2cm} & = & \hspace{-0.2cm} 4\sum_{j=1}^n\|b_{p,:}^j\|^4,\\
c_2 \hspace{-0.2cm} & = & \hspace{-0.2cm} 12{a_{\bar{p},q}^i}^T\sum_{j=1}^n\left(\|b_{p,:}^j\|^2\sum_{l=1}^{r_j}b_{p,l}^jb_{\bar{p},l}^j\right) \nonumber\\
 & = & 12{a_{:,q}^i}^T\sum_{j=1}^n\left(\|b_{p,:}^j\|^2\sum_{l=1}^{r_j}b_{p,l}^jb_{:,l}^j\right)-3c_3a_{p,q}^i, \\
 c_1 \hspace{-0.2cm} & = & \hspace{-0.2cm} 4\sum_{j=1}^n\left(\|b_{p,:}^j\|^2\left(\sum_{\substack{h=1 \\ k\neq q}}^{r_i}\sum_{l=1}^{r_j}\left({a_{:,h}^i}^Tb_{:,l}^j\right)^2+\sum_{l=1}^{r_j}\left({a_{\bar{p},q}^i}^Tb_{\bar{p},l}^j\right)^2-X_{i,j}\right)\right)+8\sum_{j=1}^n\left({a_{\bar{p},q}^i}^T\sum_{l=1}^{r_j}b_{p,l}^jb_{\bar{p},l}^j\right)^2 \nonumber\\
 \hspace{-0.2cm} & = & \hspace{-0.2cm} 4\sum_{j=1}^n\left(\|b_{p,:}^j\|^2\left(\sum_{h=1}^{r_i}\sum_{l=1}^{r_j}\left({a_{:,h}^i}^Tb_{:,l}^j\right)^2-X_{i,j}\right)\right)+8\sum_{j=1}^n\left({a_{:,q}^i}^T\sum_{l=1}^{r_j}b_{p,l}^jb_{:,l}^j\right)^2-2c_2a_{p,q}^i-3c_3{a_{p,q}^i}^2, \\
 c_0 \hspace{-0.2cm} & = & \hspace{-0.2cm} 4{a_{:,q}^i}^T\sum_{j=1}^n\left(\left(\sum_{h=1}^{r_i}\sum_{l=1}^{r_j}\left({a_{:,h}^i}^Tb_{:,l}^j\right)^2-X_{i,j}\right)\sum_{l=1}^{r_j}b_{p,l}^jb_{:,l}^j\right)-c_1a_{p,q}^i-c_2{a_{p,q}^i}^2-c_3{a_{p,q}^i}^3.
\end{eqnarray*}
\normalsize
For every entry $a_{p,q}^i$, we need to compute the different coefficients and find the root of (\ref{psd-gradient}) which minimizes the objective function (\ref{psd-contribution}).
Computing the roots of a third degree polynomial can be done in $O(1)$ operations with Cardano's method (see Appendix \ref{appA} for more details).

\subsubsection{Computational complexity of the updates}

The computation of the coefficients $c_i$'s must be implemented very carefully in order to avoid a high computational cost during the updates of the variables one after the other.
For example, we notice that the computation of the coefficient $c_0$ from scratch needs $O(nk^3)$ operations for a specific triplet $(i,p,q)$.
Updating once the $mk^2$ entries would therefore cost $O(mnk^5)$.
In the following, we explain how to reach a computational cost of $O(mk^5)$ for one pass over the $mk^2$ entries of the problem.
A loop over the `large' dimension ($n$) can be avoided with the precomputation of some quantities independent of $a_{p,q}^i$ and used during all the iterations.
For example, the term $\sum_{j=1}^n\|b_{p,:}^j\|^4$ can be precomputed and the computation of $c_3$ takes only $O(1)$ operations. 
However, the situation is more complicated for some other terms, especially $c_1$ and $c_0$. 
We describe below how to handle efficiently these computations and which quantities need to be precomputed. 

\begin{paragraph}{Computing $c_0$}
	In order to compute the coefficient $c_0$, the value of the gradient is precomputed and maintained for all the variables during the iterations.
For the purpose of clarity, we denote the quantity $\nabla_{a_{p,q}^i}f$ as $g_{p,q}^i$.
From the expression of the coefficients of (\ref{psd-gradient}), we have:
\begin{eqnarray*}
g_{p,q}^i & = & 4{a_{:,q}^i}^T\sum_{j=1}^n\left(\left(\sum_{h=1}^{r_i}\sum_{l=1}^{r_j}\left({a_{:,h}^i}^Tb_{:,l}^j\right)^2-X_{i,j}\right)\sum_{l=1}^{r_j}b_{p,l}^jb_{:,l}^j\right),\\
& = & 4{a_{:,q}^i}^T\sum_{j=1}^n\left(\left\langle A^i,B^j\right\rangle-X_{i,j}\right)B^j_{:,p} = 4{a_{:,q}^i}^TC^i_{:,p},
\end{eqnarray*}
where the different matrices $C^i=\sum_{j=1}^n\left(\left\langle A^i,B^j\right\rangle-X_{i,j}\right)B^j$, $i=1,...,m$, can be precomputed for a total of $O(mnk^2)$ operations. 
With the $C^i$ matrices available, it is possible to compute $g_{p,q}^i$ in $O(mk^3)$ operations for any triplets $(i,p,q)$.
However, $C^i$ depends on the variable $a_{p,q}^i$ and once $a_{p,q}^i$ has been assigned to its optimal value, all the entries of $C^i$ must be updated.
The entry $(u,v)$ of $C^i$ can be updated in the following way,
\footnotesize
\[
C_{u,v}^i \leftarrow C_{u,v}^i - \sum_{j=1}^n\left\langle {A^i}^{old},B^j\right\rangle B^j_{u,v} + \sum_{j=1}^n\left\langle {A^i}^{new},B^j\right\rangle B^j_{u,v} = C_{u,v}^i - \left\langle {A^i}^{old}-{A^i}^{new},\sum_{j=1}^nB^j_{u,v}B^j\right\rangle,
\]
\normalsize
and since ${A^i}^{old}-{A^i}^{new}$ is a matrix with only one non-zero row and column, the update of the $(u,v)$th entry can be done in $O(k)$ operations if the quantity $\sum_{j=1}^nB^j_{u,v}B^j$ is available.
To this end, we precompute $D_{u,v,:,:}=\sum_{j=1}^nB^j_{u,v}B^j$ for all $u$ and $v$.
To sum up, if $g^i$ is available, the coefficient $c_0$ can be computed in $O(1)$ operations.
However, after the optimization of the variable $a_{p,q}^i$, all the entries of $g^i$ and $C^i$ must be updated and it can be done in at total of $O(k^3)$ operations, which does not depend  on $n$. 

\end{paragraph}

\begin{paragraph}{Computing $c_1$} 
The second issue is the term $\sum_{j=1}^n\left({a_{:,q}^i}^T\sum_{l=1}^{r_j}b_{p,l}^jb_{:,l}^j\right)^2$ appearing in the computation of $c_1$.
The loop over the dimension $n$ can be avoided since we have
\[
E^i_{p,q}=\sum_{j=1}^n\left({a_{:,q}^i}^T\sum_{l=1}^{r_j}b_{p,l}^jb_{:,l}^j\right)^2=\left\langle a_{:,q}^i{a_{:,q}^i}^T,\sum_{j=1}^nB^j_{:,p}{B^j_{:,p}}^T\right\rangle. 
\]
In fact, if the quantity $\sum_{j=1}^nB^j_{:,p}{B^j_{:,p}}^T$ is available (and it is the case via the precomputed tensor $D$), we can maintain and update the column $q$ of $E^i$ in $O(k^2)$ operations:
\[
E^i_{l,q} \leftarrow E^i_{l,q} + 2({a^i_{p,q}}^{new}-{a^i_{p,q}}^{old}){a^i_{:,q}}^TD_{:,l,p,l}.
\]
\end{paragraph}

Table \ref{algpsd-precomputation} gathers the different quantities to precompute before the start of the iterations.
Assuming that $m$ and $n$ are of the same order of magnitude, the overall computational complexity of the precomputations is $O(mk^2\max(n,k^2))$. 
In the point of view of the space complexity, we observe that given the $A^i$'s and the $B^j$'s, the approximation matrix $\tilde{X}$ or the residual $X-\tilde{X}$ are never computed. 
In this way, the storage of a dense $m$-by-$n$ matrix is avoided (which could be impractical with a large and sparse matrix $X$). 
\small
\begin{table}
\begin{adjustbox}{center}
\small
\begin{tabular}{l|c|c}
  & computational complexity & space complexity \\
 \hline
 $A^i\leftarrow a^i{a^i}^T$ for all $i\in[m]$ & $O(mk^3)$ & $O(mk^2)$\\
 $B^j\leftarrow b^j{b^j}^T$ for all $j\in[n]$ & $O(nk^3)$ & $O(nk^2)$\\
 $C^i\leftarrow\sum_{j=1}^n\left(\left\langle A^i,B^j\right\rangle-X_{i,j}\right)B^j$ for all $i\in[m]$, & $O(mnk^2)$ & $O(mk^2)$\\
 $g^i \leftarrow 4{a^i}^TC^i$ for all $i\in[m]$ & $O(mk^3)$ & $O(mk^2)$\\
 $D_{u,v,:,:} \leftarrow \sum_{j=1}^nB^j_{u,v}B^j$ for all $u,v=1,...,k$ & $O(nk^4)$ & $O(k^4)$\\
$E^{i}_{p,q} \leftarrow \left\langle a_{:,q}^i{a_{:,q}^i}^T,\sum_{j=1}^nB^j_{:,p}{B^j_{:,p}}^T\right\rangle$ for all $i\in[m]$, $p,q\in[k]$ & $O(mk^4)$ & $O(mk^2)$
\end{tabular}
\end{adjustbox}
\caption{List of precomputations for the CD methods}
\label{algpsd-precomputation}
\end{table}
\normalsize

%

\subsubsection{Variables selection: cyclic or greedy}
Algorithm~\ref{algcyclicCD} illustrates a cyclic run of a CD scheme over all the variables.
After the computation of the optimal value of one of the $mk^2$ entries of the problem, the updates of $C^i$ and $g^i$ in $O(k^3)$ operations are the bottleneck of the method causing the overall $O(mk^5)$ complexity.

\begin{algorithm}[h!]
\caption{Optimize subproblem~\eqref{psd-obj-newvariables} (cyclic coordinate descent)}
\label{algcyclicCD}
\begin{algorithmic}[1]
\STATE INPUT: $X \in \mathbb{R}_+^{m\times n}$, $\{b^1,...,b^n\} \in \mathbb{R}^{k\times r}$.
\STATE OUTPUT: $\{a^1,...,a^m\} \in \mathbb{R}^{k\times r}$.
\STATE $\{a^1,...,a^m\} \leftarrow$ Initialization
\STATE $[C,D,E,g]\leftarrow$Precomputation$(X,\{a^1,...,a^m\},\{b^1,...,b^n\})$
\FOR{$i=1:m$}
\FOR{$p=1:k$}
\FOR{$q=1:r$}
\STATE $x\leftarrow a^i_{p,q}$
\STATE $c_3\leftarrow 4D_{p,p,p,p}$
\STATE $c_2\leftarrow 12{a^i_{:,q}}^TD_{p,p,p,:}-3c_3x$
\STATE $c_1\leftarrow 4C^i_{p,p}+8E^i_{p,q}-2c_2x-3c_3x^2$
\STATE $c_0\leftarrow 4g^i_{p,q}-c_1x-c_2x^2-c_3x^3$
\STATE ${a_{p,q}^i} \leftarrow CardanoMethod(c_3,c_2,c_1,c_0)$ \label{psd-cyclic-line13}
\STATE Update $C^i$, $g^i$, and $E^i_{:,q}$
\ENDFOR
\ENDFOR
\ENDFOR
\end{algorithmic}
\end{algorithm}

As explained above, the gradient of any variable is always available in Algorithm~\ref{algcyclicCD}.  
In order to improve the efficiency of the algorithm, we propose to use the information given by the gradient for selecting first the coordinates in a greedy way instead of processing them cyclically.
This is called the Gauss-Southwell rule: at each iteration, the variable with the largest gradient is updated. It allows to guide the CD scheme towards the coordinates that will potentially decrease the objective function the most. 
Algorithm~\ref{algGSCD} describes the implementation of the Gauss-Southwell strategy for PSD factorization.
The main difference with Algorithm~\ref{algcyclicCD} lies in the selection of the variables to optimize.

\begin{algorithm}[h!]
\caption{optimizesubproblem (Gauss-Southwell coordinate descent)}
\label{algGSCD}
\begin{algorithmic}[1]
\STATE INPUT: $X \in \mathbb{R}_+^{m\times n}$, $\{b^1,...,b^n\} \in \mathbb{R}^{k\times r}$, $\alpha \in \mathbb{R}_+$.
\STATE OUTPUT: $\{a^1,...,a^m\} \in \mathbb{R}^{k\times r}$.
\STATE $\{a^1,...,a^m\} \leftarrow$ Initialization
\STATE $[C,D,E,g]\leftarrow$Precomputation$(X,\{a^1,...,a^m\},\{b^1,...,b^n\})$
\FOR{$i=1:m$}
\FOR{$t=1:\lceil\alpha kr\rceil$}
\STATE $(p^{\ast},q^{\ast})=\arg \max_{p,q} |g^i_{p,q}|$
\STATE $x\leftarrow a^i_{p{\ast},q{\ast}}$
\STATE $c_3\leftarrow 4D_{p{\ast},p{\ast},p{\ast},p{\ast}}$
\STATE $c_2\leftarrow 12{a^i_{:,q{\ast}}}^TD_{p{\ast},p{\ast},p{\ast},:}-3c_3x$
\STATE $c_1\leftarrow 4C^i_{p{\ast},p{\ast}}+8E^i_{p{\ast},q{\ast}}-2c_2x-3c_3x^2$
\STATE $c_0\leftarrow 4g^i_{p{\ast},q{\ast}}-c_1x-c_2x^2-c_3x^3$
\STATE ${a_{p^{\ast},q^{\ast}}^i} \leftarrow CardanoMethod(c_3,c_2,c_1,c_0)$
\STATE Update $C^i$, $g^i$, and $E^i_{:,q}$
\ENDFOR
\ENDFOR
\end{algorithmic}
\end{algorithm}

We propose to make a number of iterations on each factor $i$ proportional to $kr$ (the number of variables) using the parameter $\alpha$. 
In Section~\ref{sec-results}, the performances of Algorithm~\ref{algGSCD} are compared for different values of $\alpha$.

\subsubsection{Inner rank of the factors} \label{innerrank}

In many cases, the factors  $A^i$'s and the $B^j$'s are rank deficient.
For example, in the exact case ($X_{ij} = \langle A^i,B^j\rangle$ for all $i,j$), 
if $X_{ij} = 0$ and the $i$th row of $X$ and $j$th column of $X$ are not identically zero (implying $A^i \neq 0$ and $B^j \neq 0$), $A^i$ and $B^j$ cannot be full rank otherwise $\langle A^i,B^j\rangle > 0$.
For slack matrices, there is at least one zero per row and per column in $X$, hence $r_i \leq k-1$ for all $i$. 
In fact, this idea can be generalized~\cite{lee2012support} to improve the upper bound on the $r_i$'s, and was used for example in~\cite{fawzi2016rational}. 

With the CD methods previously presented, it is easy to allow different values for the rank of the $A^i$'s by using initial factors $a^i$'s with appropriate sizes.
However, for the numerical experiments in Section \ref{sec-results}, we will use $r_i = k$ for all $i$ to have a fair comparison with FPGM and to check whether the coordinate descent algorithms are able to generate low-rank factors. 
Moreover, this possibility to handle rank deficient factors will allow us to focus on the problem of the square root rank where $r_i = 1$ for all factors; see Section \ref{square-root-rank}.

\section{Numerical experiments} \label{sec-results}
The algorithms presented in the previous section are the first numerical methods developed for solving the optimization problem (\ref{psd-opt}).
It is therefore not possible to any make experimental comparisons with algorithms from the literature.
However, this section has two main goals:
\begin{itemize}

	\item In Algorithms \ref{algfgradient} and \ref{algGSCD}, there are parameters that may influence the effectiveness of the methods, $\Delta$ and $\alpha$ respectively.
    Hence the first goal is to compare the performances of these two algorithms for different values of the parameters. 
		
    \item Once the best values of the parameters are known, the second goal is to compare the Algorithms \ref{algfgradient}, \ref{algcyclicCD} and \ref{algGSCD}.    
	This will allow us to select the most effective algorithm to solve the PSD factorization problems discussed in Section~\ref{sec-application}.
		
\end{itemize}

\subsection{Initialization and scaling} \label{initscaling}

Algorithms \ref{algfgradient}, \ref{algcyclicCD} and \ref{algGSCD} are iterative and need starting points.
In this paper, the entries of the $a^i$'s and the $b^j$'s are initialized using the normal distribution $\mathcal{N}(0,1)$.
Note that for Algorithm~\ref{algfgradient}, we use $a^i{a^i}^T$ and $b^j{b^j}^T$ as random initial iterates so that all algorithms are initialized with the same values.  

However, it may happen that with such random factors, we have an initial approximation matrix $\tilde{X}$ way larger or smaller than $X$.
In order to avoid such situations, we scale the initial factors compared to $X$: given intial iterates $a^i$ and $b^j$, we compute
\[
\lambda^{\ast}=\arg\min_{\lambda}\sum_{i=1}^m\sum_{j=1}^n\left(X_{i,j}-\lambda\left\langle A^i,B^j\right\rangle\right)^2 =\arg\min_{\lambda} \left\|X - \lambda\mathcal{A}^T\mathcal{B}\right\|^2_F=\frac{\left\langle X\mathcal{B}^T,\mathcal{A}\right\rangle}{\left\langle\mathcal{B}\mathcal{B}^T,\mathcal{A}\mathcal{A}^T\right\rangle},
\]
with $\mathcal{A}=\begin{pmatrix}vec(A^1) \hspace{0.1cm} ... \hspace{0.1cm} vec(A^m)\end{pmatrix}$ and $\mathcal{B}=\begin{pmatrix}vec(B^1) \hspace{0.1cm} ... \hspace{0.1cm} vec(B^n)\end{pmatrix}$.
The initial error is therefore
\[
	e_0=\left\|X - \tilde{X}\right\|_F=\sqrt{\|X\|^2_2-\frac{\left\langle X\mathcal{B}^T,\mathcal{A}\right\rangle^2}{\left\langle\mathcal{B}\mathcal{B}^T,\mathcal{A}\mathcal{A}^T\right\rangle}} \leq \left\|X\right\|_F,
\]
with the appropriate scaling,
\begin{itemize}
	\item $A^i\leftarrow\lambda^{\ast}A^i$ for $i=1,...,m$ for FPGM and,
	\item $a^i\leftarrow\sqrt{\lambda^{\ast}}a^i$ for $i=1,...,m$ for the CD methods.
\end{itemize}

\subsection{Data sets} \label{psd-datasets}

The matrices used for the numerical comparisons are slack matrices; see the discussion in Section~\ref{sec-defmot}.   
Table~\ref{psd-datasets-table} summarizes the different matrices used in the tests.
The factorization rank $k$ used in the experiments is specified in the fourth column.
Note that this is not necessarily the true value of the $\rankp$ which is used, but it is either a conjecture or an upper bound.
The data set is composed of three types of matrices: 
\begin{itemize}

	\item The slack matrices of regular $n$-gons are $n$-by-$n$ circulant matrices for which the $(i,j)$th entry is the slack between the $i$th facet and the $j$th vertex of the regular $n$-gon (see \cite{VGG15} for more details on the construction of such matrices).
	The values of the factorization rank $k$ are given by the conjecture made on the $\rankp$ 
	of regular $n$-gons in Section~\ref{conjecturepsdrank}. 
	
    \item For a given positive integer $n$, let $U_n$ (resp.\@ $V_n$) be the $\{0,1\}^{{{n}\choose{\lfloor\frac{n}{2}\rfloor}}\times n}$ (resp.\@ $\{0,1\}^{{{n}\choose{\lceil\frac{n}{2}\rceil}}\times n}$) matrix where the rows correspond to the subsets of $\{1,...,n\}$ of size ${n}\choose{\lfloor\frac{n}{2}\rfloor}$ (resp.\@ ${n}\choose{\lceil\frac{n}{2}\rceil}$).
    Let $P_n$ be the ${n}\choose{\lfloor\frac{n}{2}\rfloor}$-by-${n}\choose{\lceil\frac{n}{2}\rceil}$ matrix defined as
    \[
		P_n=U_nV_n^T.
    \]
    These matrices have an interpretation in terms of an inscribed polytope in the $(n-2)$-sphere (see Problems 9.1 and 9.2 in \cite{FGP14}).
    The exact value of $\rankp(P_n)$ is not known but it is bounded as follows,
    \[
    	\left\lceil\frac{\sqrt{1+8n}-1}{2}\right\rceil\leq \rankp(P_n)\leq 2\left\lceil \sqrt{n} \, \right\rceil,
    \]
    except for $n=5$ for which $3\leq \rankp(P_5)\leq 4$.
    The values of the factorization rank $k$ of the matrices $P_n$ used in the tests are the upper bounds mentioned above.
    \item 
The correlation polytope is the convex hull of all $n$-by-$n$ rank-one 0/1 matrices.
Let $COR_n$ be a submatrix of the slack matrix of the correlation polytope.
The rows and columns of this $2^n$-by-$2^n$ matrix are indexed by vectors $u, v \in \{0,1\}^n$ such that 
\[
COR_n(u,v) 
= \left(1-u^Tv\right)^2.
\]
Although the nonnegative rank of $COR_n$ has been proved to be exponential in $n$, there exists an explicit PSD factorization such that $\rankp(COR_n)=n+1$~\cite{FMPTdW12}. 
These values are used for the factorization rank in the tests.
\end{itemize}

\begin{table}[ht]
\begin{center}
\begin{tabular}{|c||c|c|c|c|}
\hline
& $m$  &   $n$ & $k$  \\ 
\hline  
slack matrix of the $12$-gon  &	12 & 12  & 5  \\ 
slack matrix of the $16$-gon  &	16 & 16  & 5  \\ 
slack matrix of the $20$-gon  &	20 & 20  & 6  \\ 
slack matrix of the $24$-gon  &	24 & 24 & 6  \\ 
slack matrix of the $28$-gon  &	28 & 28 & 6 \\ 
slack matrix of the $32$-gon  &	32 & 32 & 6\\
\hline  
$P_5$ &	 10   & 10  & 4 \\ 
$P_6$ &	 20   & 20  & 6 \\ 
$P_7$ &	 35   & 35  & 6 \\ 
\hline 
$COR_3$ &	 8   & 8  & 4   \\ 
$COR_4$ &	 16   & 16  & 5 \\ 
$COR_5$ &	 32   & 32  & 6 \\ 
\hline 
\end{tabular}
\caption{Benchmark of nonnegative matrices used in the numerical comparisons.} 
\label{psd-datasets-table}
\end{center} 
\end{table}


\subsection{Comparisons for different values of the parameters}


In order to compare the performances of the algorithms, we use the measure $E(t)$ defined by
\begin{equation} \label{psd-average}
E(t)=\frac{e(t)}{e_0} 
\end{equation}
where $e_0$ is the initial error (see Section \ref{initscaling}),
and $e(t)$ is the error $\|X-\tilde{X}\|_F$ achieved by an algorithm for a given initialization within $t$ seconds. Since our algorithms are nonincreasing, we have $E(t) \in [0,1]$ for all $t$, with $E(0) = 1$ and $E(t) \rightarrow_{t \rightarrow \infty} 0$ if the corresponding algorithm converges towards an exact factorization. 
In order to illustrate the efficiency of a given algorithm, (\ref{psd-average}) has the advantage that it makes sense to take the average of $E(t)$ for several initializations and data sets and display a single curve.
The algorithms were run 10 times with different initializations during 60 seconds for the following parameters values:
\begin{itemize} 
	\item $\Delta=\{\frac{1}{k},0.5,1,5,10,20,30\}$ for Algorithm~\ref{algfgradient}, and
	\item $\alpha=\{0.05,0.1,0.5,1,5,20,40\}$ for Algorithm~\ref{algGSCD}.
\end{itemize}
FPGM was implemented with Matlab while the CD methods were developed in C with a Matlab interface using Mex files.
The reason is that Matlab is not a well-suited language when one requires to perform many loops as in Algorithms \ref{algcyclicCD} and \ref{algGSCD}.
The codes are available at \url{https://sites.google.com/site/exactnmf/}.  
All tests were performed on a PC Intel CORE i5-4570 CPU @3.2GHz $\times$ 4, with 7.7G RAM.

\begin{figure*}
	\begin{center}
    \hspace*{-0.4cm} 
	\begin{tabular}{cc}
	\includegraphics[width=8.5cm]{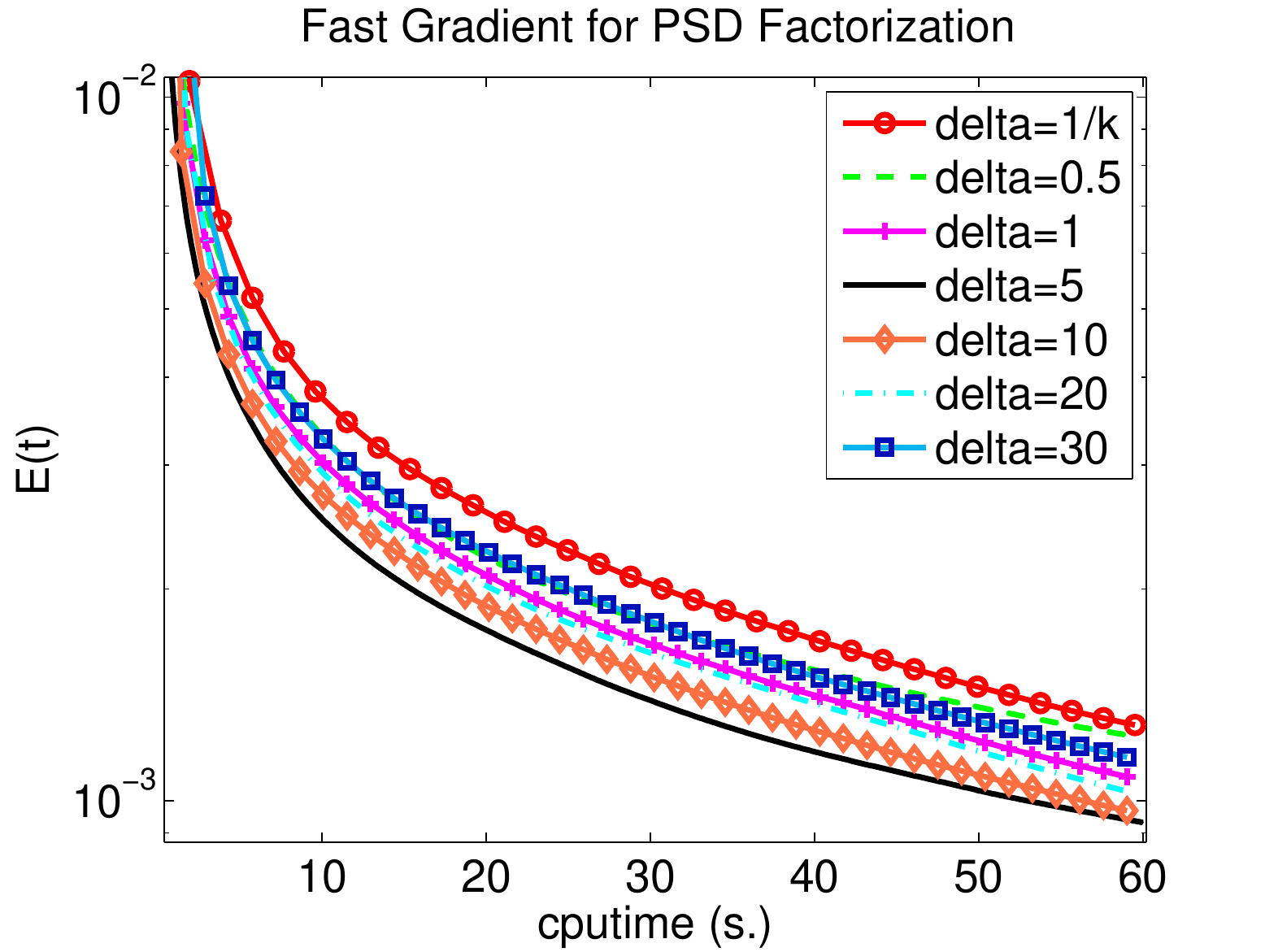} & \includegraphics[width=8.5cm]{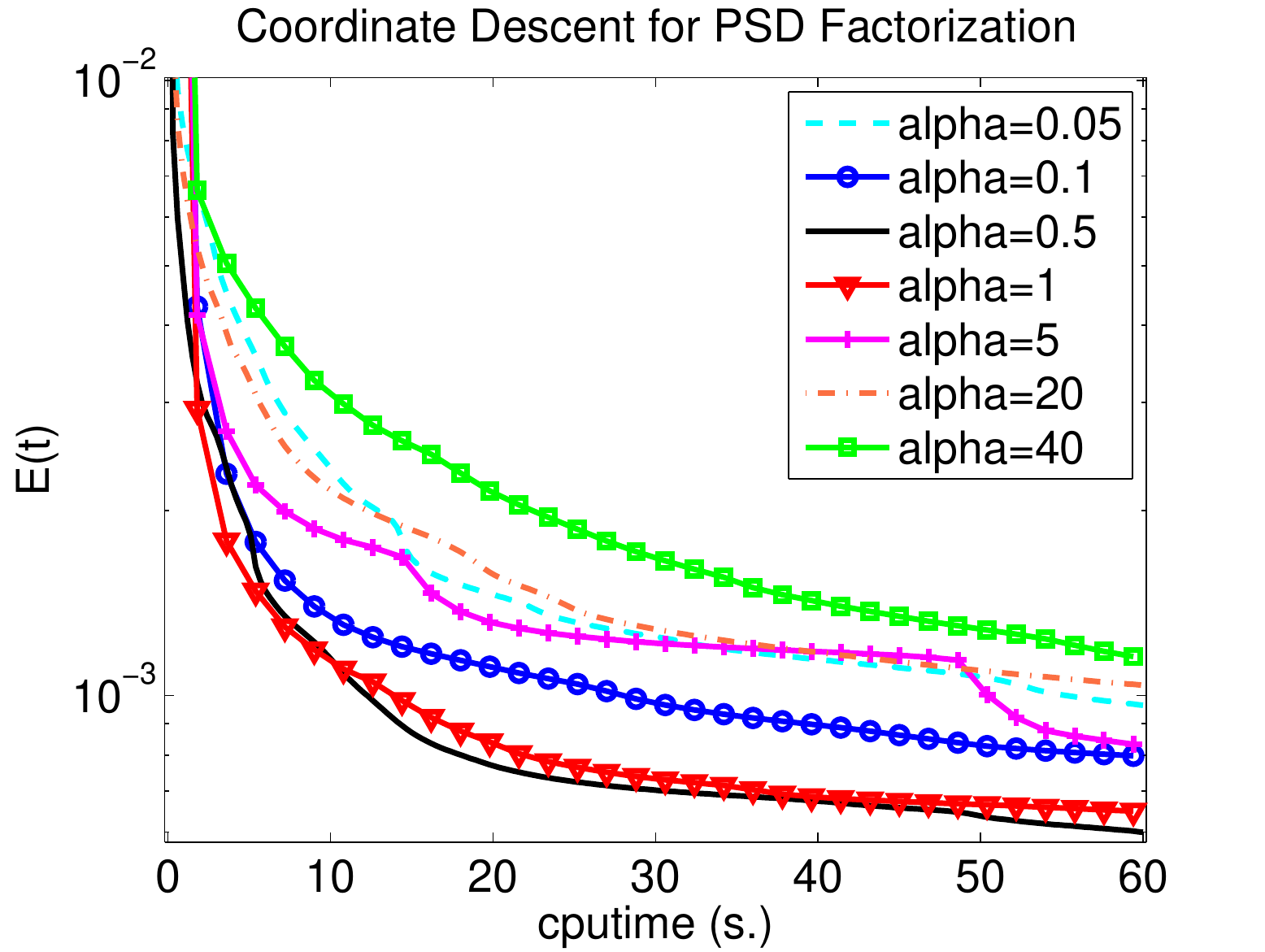}\\
    (a) & (b)
	\end{tabular} 
\end{center}
\caption{Evolution of the average measure $E(t)$ for different values of the parameters $\Delta$ and $\alpha$ on the data sets of Table \ref{psd-datasets-table}.} 
\label{psd-comp-algo} 
\end{figure*}

The results are displayed on Figure \ref{psd-comp-algo}. 

For FPGM, we observe that the number of inner steps does not influence the efficiency significantly. 
We observe that the best average performances are obtained around $\Delta=5$.
For the Gauss-Southwell algorithm, the best value of the parameter $\alpha$ is between $0.5$ and $1$.
It means that the number of updated entries must be roughly the same as in the cyclic case. 
For the numerical tests that follow, we use the following algorithms:
\begin{itemize}
	\item FPGM with $\Delta=5$.
    \item The cyclic CD algorithm.
    \item The Gauss Southwell CD algorithm with $\alpha=0.5$.
\end{itemize}

In the remaining of the section, we compare the performances of these algorithms.
Instead of ploting an average measure, for each matrix and each method, we display the curves of the error $\|X-\tilde{X}\|_F$ corresponding to five different initializations. 
It allows us to observe the behavior of the methods for different starting points.
The data sets used are those described in Table \ref{psd-datasets-table}.
For each type of matrices, we present the results for two instances: the matrices with the smallest and the largest size.

\begin{figure*}
	\begin{center}
    \hspace*{-0.4cm} 
	\begin{tabular}{cc}
    \includegraphics[width=8.5cm]{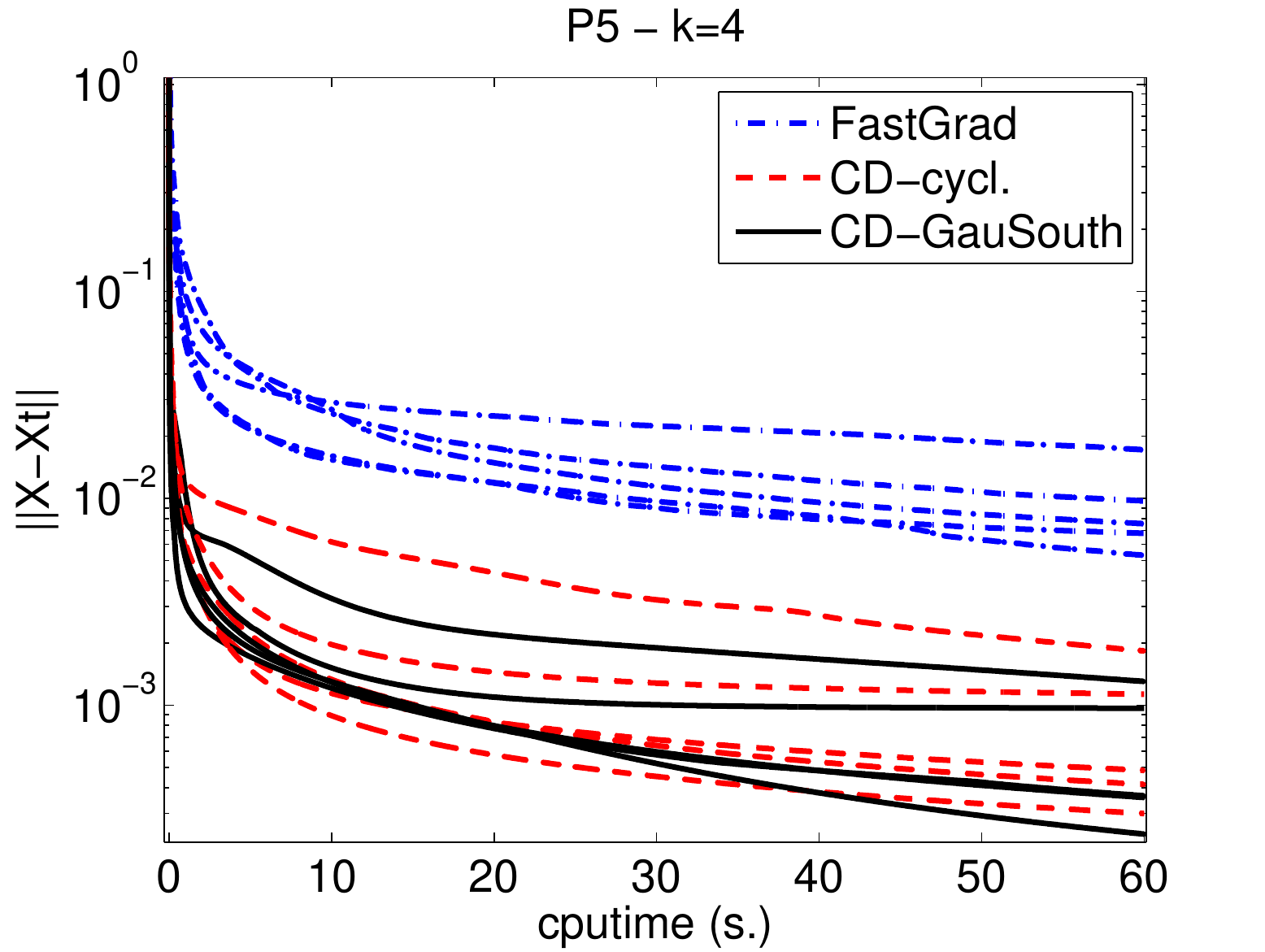} &	\includegraphics[width=8.5cm]{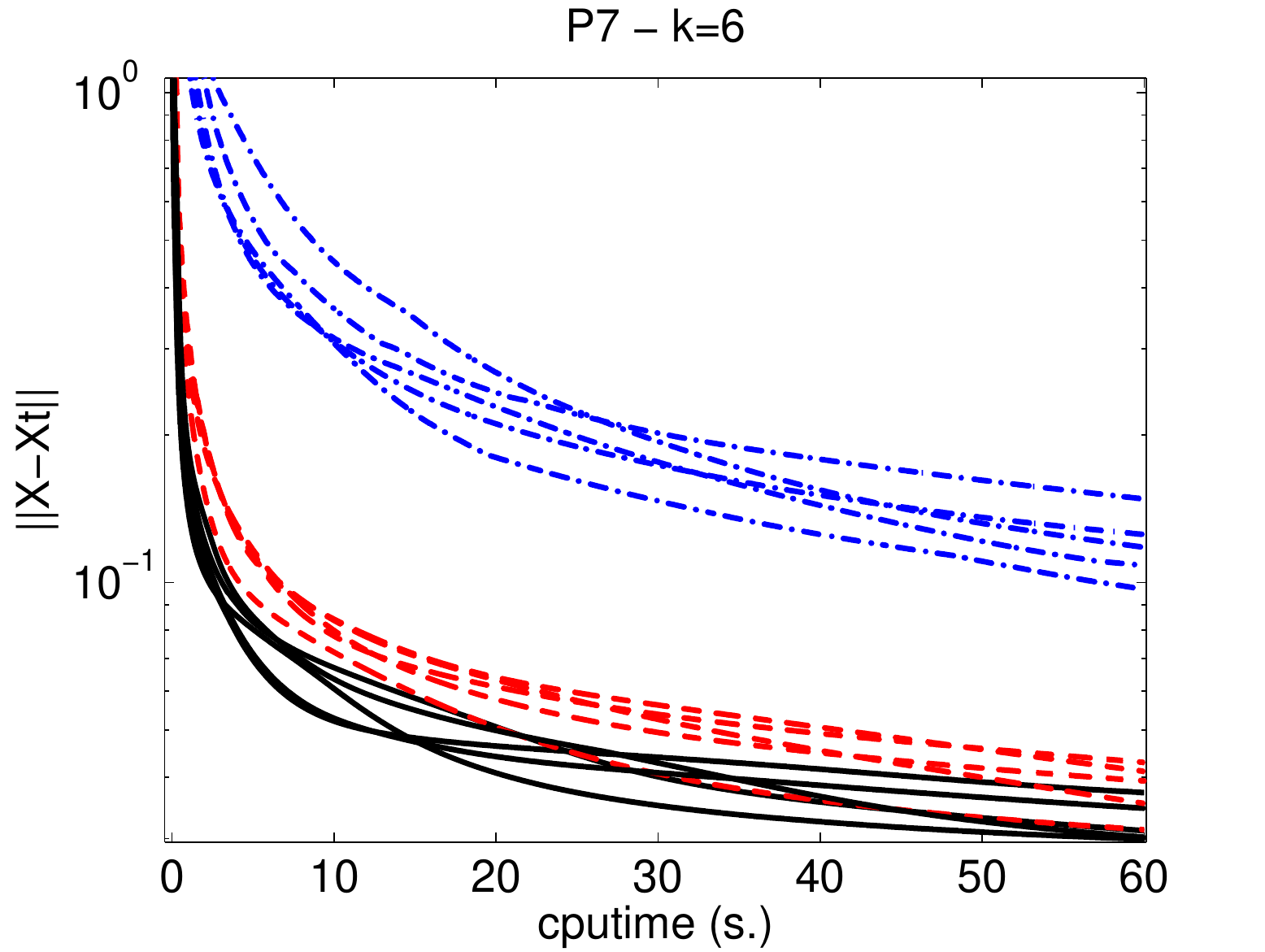}  \\
	\includegraphics[width=8.5cm]{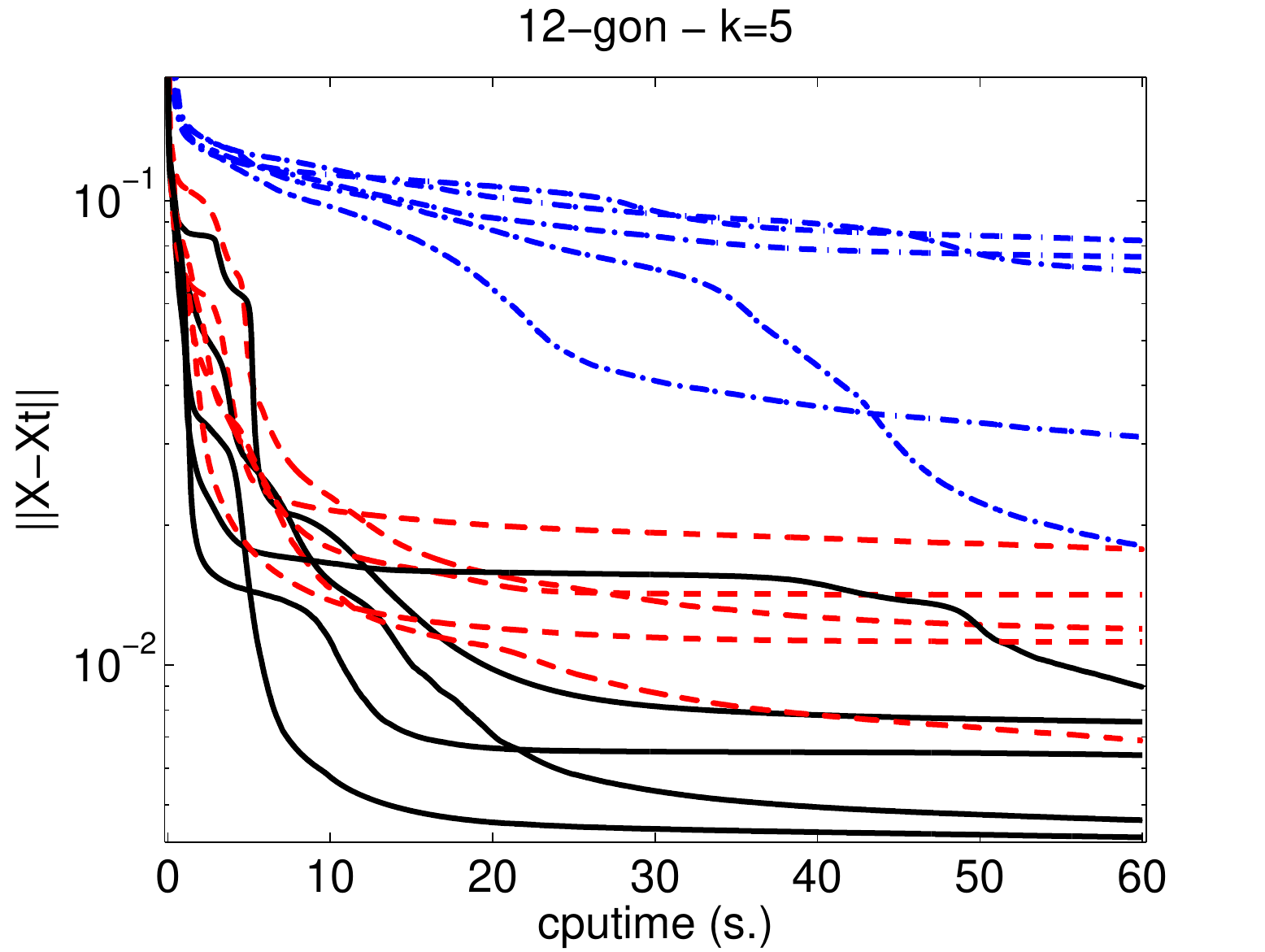} & \includegraphics[width=8.5cm]{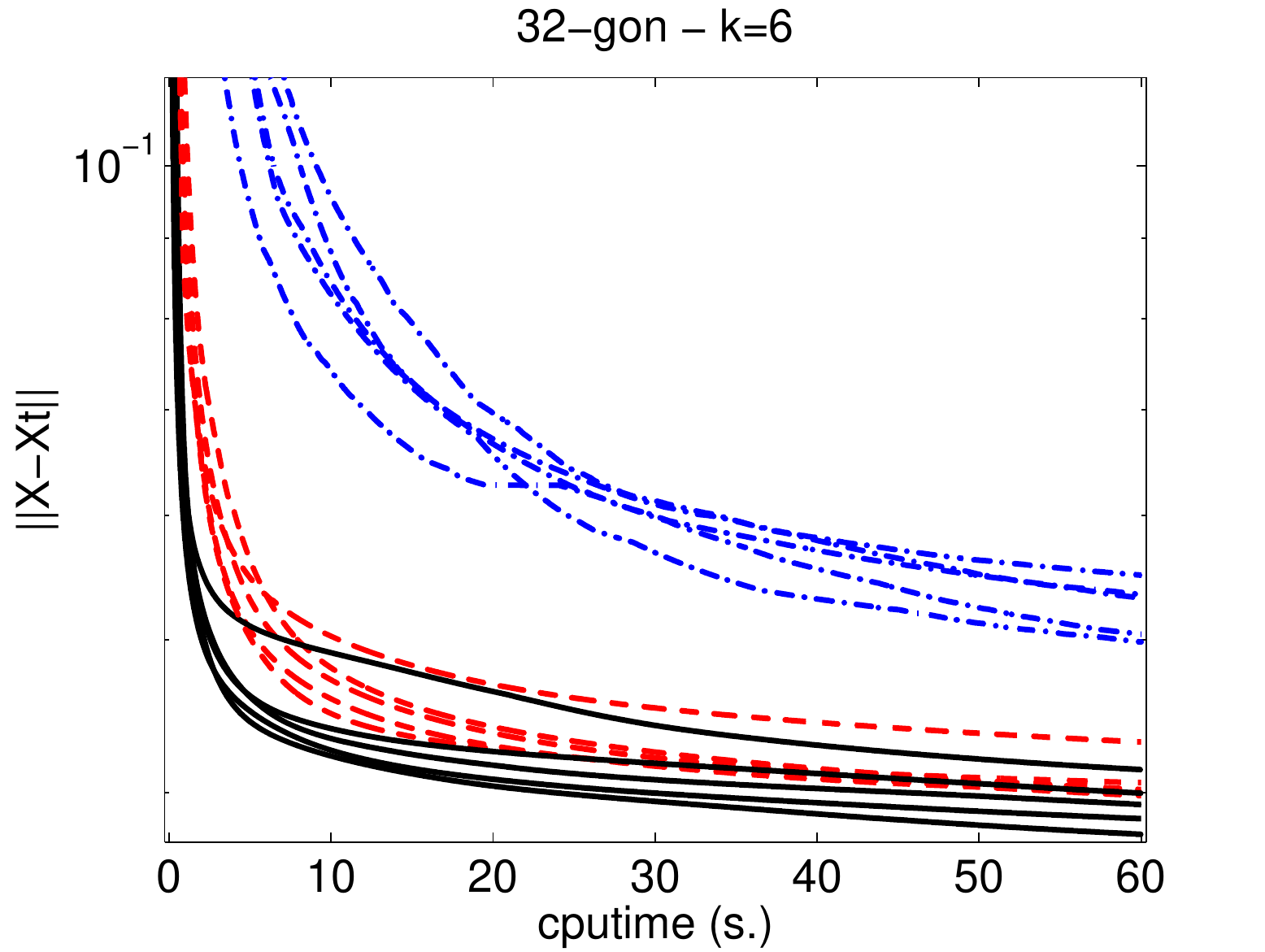} \\
	  \includegraphics[width=8.5cm]{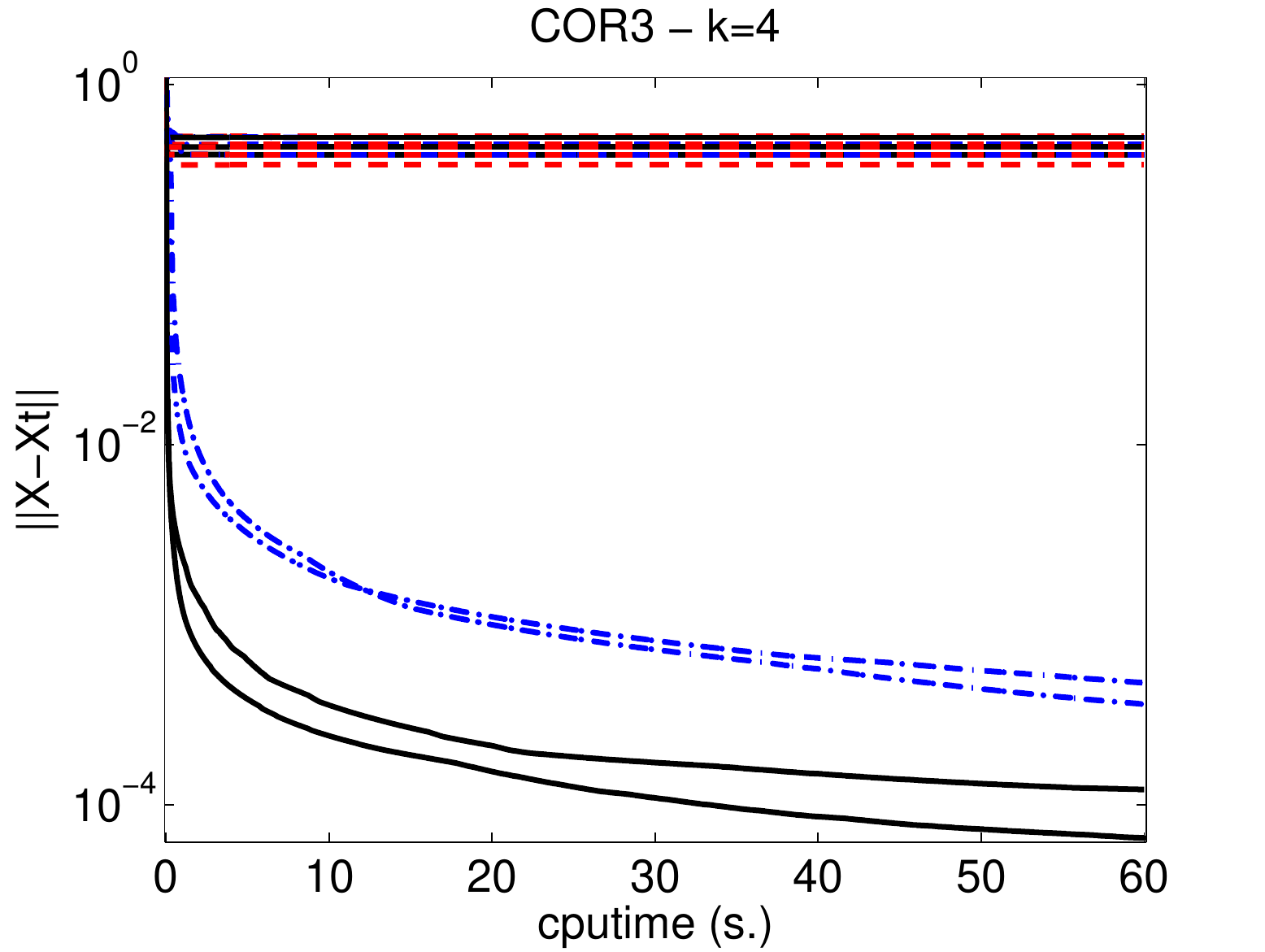} & \includegraphics[width=8.5cm]{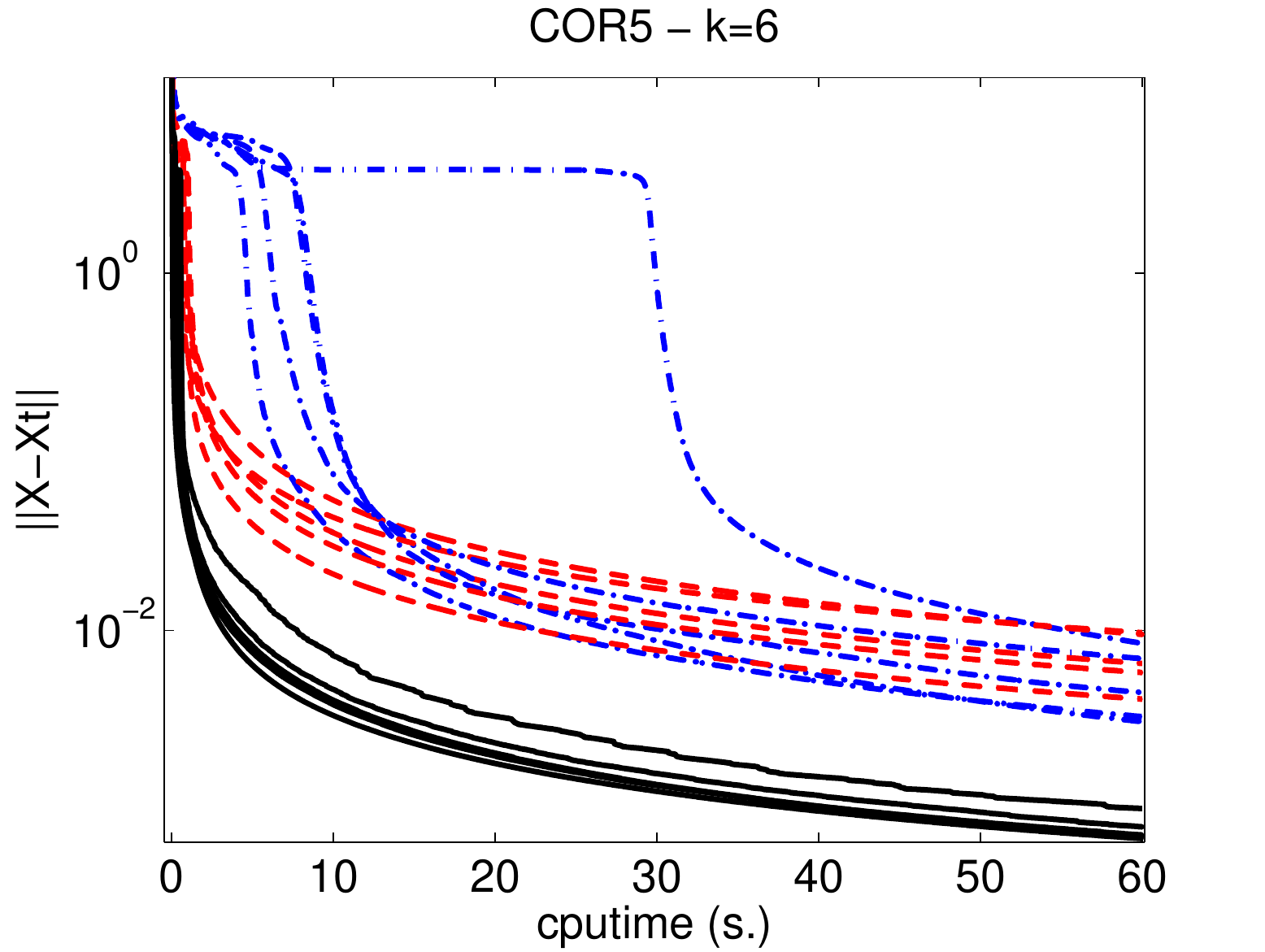} 
	\end{tabular} 
\end{center}
\caption{Evolution of the error for the different algorithms on the dataset.} 
\label{psd-comp-dataset} 
\end{figure*}

From Figure \ref{psd-comp-dataset}, we observe the following:
\begin{itemize}
	\item There is a general trend emerging from these numerical tests: the Gauss-Southwell CD method outperforms the cyclic strategy, while this last method performs better than FPGM.
    \item Algorithm are very sensitive to initialization. 
    For example, the solutions obtained with FPGM on the $12$-gon after 60 seconds are rather different, illustrating the fact the local algorithms can get stuck in local minima. 
	This is clear from the results obtained with the $COR3$ matrix where most of the runs get stuck in local minima. 
    \item Although the results presented on the left of Figure~\ref{psd-comp-dataset} are instances of small sizes, the final error $\|X-\tilde{X}\|_F$ remains relatively large even after 60 seconds.
    It contrasts with NMF where the convergence on small matrices is faster~\cite{VGGT14}.   
\end{itemize} 

In conclusion, we recommend to use the Gauss-Southwell CD method which performs best in most cases.
This algorithm will therefore be used in the next section for several applications where the $\rankp$ is sought.

\section{Applications} \label{sec-application}

In this section, we discuss the use of our numerical algorithms for the computation of the psd-rank of particular matrices.
In this purpose, let us give the following (obvious) fact.
\begin{observation}
For a given matrix $X$ with $\rankp(X)=k$, let us denote $\tilde{X}^{\ast}_l$ the best approximation matrix with $l$-by-$l$ PSD factors.
By definition of $\rankp$, we have
\begin{eqnarray*}
	\|X-\tilde{X}^{\ast}_l\|_F=0 \text{ for all }l\geq k, \quad  \text{ and } \quad 
	\|X-\tilde{X}^{\ast}_l\|_F>0 \text{ for all }l< k.
\end{eqnarray*}

\end{observation} 
Given a matrix $X$ and a target factorization rank $k$, our nonlinear local optimization methods provide no guarantee; we can only hope to identify good local minima of the nonconvex problem (\ref{psd-opt}).
However, as experimentally demonstrated in \cite{VGGT14} for exact NMF, such algorithms can be used in multi-start strategies to detect if the error $\|X-\tilde{X}\|_F$ gets (close) to zero. 
Moreover, beside conjectures on the psd-rank, Algorithms \ref{algcyclicCD} and \ref{algGSCD} can be helpful to find exact factorizations by trial and error and by fixing manually some entries to specific values. 

As an illustration, we discuss the value of the psd-rank of the regular polygons in Section \ref{conjecturepsdrank}.
With the help of Algorithm~\ref{algGSCD}, a conjecture is proposed which is confirmed showing exact factorizations, for the first time, for $n=5$, $n=8$ and $n=10$. 
In Sections \ref{completelypsd} and \ref{square-root-rank}, we show how to adapt our methods in order to deal with two related problems, the completely PSD Factorization problem and the problem of computing the square root rank.

%
%

\subsection{Conjecture on the psd-rank of regular $n$-gons} \label{conjecturepsdrank}

Let $S_n$ denote the slack matrix of the regular $n$-gon.
We have that
\begin{equation*}\label{lb-ub-ngon}
\Omega\left(\frac{\log n}{\log\log n}\right)\leq \rankp(S_n)\leq 2\lceil \log_2(n) \rceil,
\end{equation*}
where the first inequality comes from quantifier elimination theory~\cite{GPT13, GRT13} and the second inequality uses the upper bound on $\rank_+(S_n)$~\cite{FRT12}. 
The exact value of $\rankp(S_n)$ is unknown for general $n$. However, it is known that 
(i)~the psd-rank of the square is three, 
(ii)~all pentagons and hexagons have psd-rank exactly four 
and 
(iii)~the psd-rank of the heptagons is either four or five~\cite{GRT13}. 
Moreover, to the best of our knowledge, an explicit factorization for regular $n$-gons is only known for $n=3$, $n=4$ and $n=6$.

For different values of $n$ and $k$, we run Algorithm~\ref{algGSCD} on $S_n$ with the inner rank of the factors $r=k-2$ (see Section \ref{innerrank}). Table \ref{psd-ngons} reports the smallest relative error $\frac{\|X-\tilde{X}\|_F}{\|X\|_F}$ found after 100 runs of 10 seconds with different initializations. 
\begin{table}[ht] 
\begin{center}
\begin{tabular}{|c|||c|c|c|c|c|}
\hline
       &  $k=3$            & $k=4$            & $k=5$           & $k=6$           & $k=7$  \\ 
\hline
$n=3$  &   \textbf{3.1e-7} &                  &                 &                 &        \\
$n=4$  &   \textbf{1.3e-7} & 6.6e-7           &                 &                 &        \\
$n=5$  &   0.065           & \textbf{2.4e-6}  & 3.8e-6          &                 &        \\
$n=6$  &   0.049           & \textbf{5.5e-6}  & 5.4e-6          &                 &        \\  
$n=7$  &   0.036           & \textbf{3.9e-5}  & 1.3e-5          &                 &        \\  
$n=8$  &   0.028           & \textbf{1.9e-5}  & 3.4e-5          &                 &        \\  
$n=9$  &   0.022           &  0.004           & \textbf{8.5e-5} & 3.7e-5          &        \\  
$n=10$ &   0.018           &  0.003           & \textbf{8.1e-5} & 4.9e-5          &        \\  
$n=11$ &   0.015           &  0.006           & \textbf{1.4e-4} & 5.4e-5          &        \\  
$n=12$ &   0.012           &  0.007           & \textbf{2.7e-4} & 1.2e-4          &        \\  
$n=13$ &   0.01            &  0.007           & \textbf{5.5e-4} & 1.6e-4          &        \\
$n=14$ &   0.009           &  0.006           & \textbf{6.9e-4} & 2.8e-4          &        \\
$n=15$ &   0.008           &  0.005           & \textbf{8e-4}   & 4e-4            &        \\
$n=16$ &   0.007           &  0.005           & \textbf{0.001}  & 5e-4            &        \\
$n=17$ &   0.006           &  0.004           & 0.002           & \textbf{5.6e-4} & 4.6e-4 \\
\hline 
\end{tabular}
\caption{Smallest relative error obtained over $100$ different runs of $10$ seconds.}
\label{psd-ngons}
\end{center} 
\end{table}

In order to guess a value for the psd-rank of $S_n$, we have to look at the corresponding row of Table \ref{psd-ngons}.
If an exact factorization is possible for $k$, the error should be close to zero in the entry $(n,k)$ and larger in the entry $(n,k-1)$. 
For the smallest regular $n$-gons ($n=3,4,5,$ and $6$), the obtained errors are consistent with the known values of the psd-rank.
For $n=7$ and $n=8$, the results 
suggest\footnote{Example~\ref{ex8gon} provides an explicit PSD factorization of size 4 for $S_8$. For $S_7$, we were not able to obtain such an exact factorization of size 4, although we have tried many different initializations. 
It is possible that $\rankp(S_7) = 5$ since there is no result about the monotonicity of the PSD rank of regular $n$-gons (this is, as far as we know, an open question). In fact, \cite{goucha2016ranks} showed that monotonicity does not hold for the PSD rank over the complex numbers with $\rank_{\psd}^{\mathbb{C}}(S_6) = 3 < 4 \leq \rank_{\psd}^{\mathbb{C}}(S_5)$.} 
that $\rankp(S_7)=\rankp(S_8)=4$. 
Actually, there is a pattern emerging for $n\geq 7$ leading to the following conjecture:
\begin{conj} \label{conj-psd-ngons}
The psd-rank of $S_n$, the slack matrix of the regular $n$-gon, is given by 
$$\rankp(S_n)=1 + \lceil \log_2(n)\rceil.$$
\end{conj} 

In Table \ref{psd-ngons}, the entries corresponding to the conjecture are highlighted in bold.
We have not pursued the computations beyond $n>17$ because the results are less and less clear.
The reason is that as $n$ gets bigger, the regular $n$-gon get closer to the circle which has a psd-lift of size $2$.

By trial and error and by fixing more and more entries manually in the factors, we were able to construct, for the first time, an exact PSD factorization of the $5$-gon, the $8$-gon and the $10$-gon with respective sizes consistent with Conjecture~\ref{conj-psd-ngons}; see the examples below. 
\begin{exmp}
With $\phi=\frac{1+\sqrt{5}}{2}$, a slack matrix of the regular $5$-gon is given by:
$$S_5=\begin{pmatrix}
0 & 1 & \phi & 1 & 0 \\
0 & 0 & 1 & \phi & 1 \\
1 & 0 & 0 & 1 & \phi \\
\phi & 1 & 0 & 0 & 1 \\
1 & \phi & 1 & 0 & 0
\end{pmatrix}.$$ 
A $\mathcal{S}^4_+$-factorization of $S_5$ is given by the following factors:
\small
$$a^i=\left\{\begin{pmatrix}
1 \\
0 \\
0 \\
0
\end{pmatrix},
\begin{pmatrix}
0 \\
1 \\
0 \\
0
\end{pmatrix},
\begin{pmatrix}
0 \\
0 \\
1 \\
0
\end{pmatrix},
\begin{pmatrix}
0 \\
0 \\
0 \\
1
\end{pmatrix},
\begin{pmatrix}
\frac{1}{1+\sqrt{\phi}} \\
\frac{1}{1+\sqrt{\phi}} \\
1 \\
\phi-\frac{1}{\sqrt{\phi}}
\end{pmatrix}\right\},$$
$$
b^j=\left\{
\begin{pmatrix}
0 & 0 \\
0 & 0 \\
\frac{1-(\sqrt{\phi})^3}{2} & \sqrt{1-\left(\frac{1-(\sqrt{\phi})^3}{2}\right)^2}\\
\sqrt{\phi} & 0
\end{pmatrix},
\begin{pmatrix}
1 \\
0 \\
0 \\
1
\end{pmatrix},
\begin{pmatrix}
\sqrt{\phi} \\
1 \\
0 \\
0
\end{pmatrix},
\begin{pmatrix}
1 \\
a \\
-1 \\
0
\end{pmatrix},
\begin{pmatrix}
0 \\
-1 \\
\sqrt{\phi} \\
-1
\end{pmatrix}\right\}.
$$

\normalsize
\end{exmp}

\begin{exmp} \label{ex8gon} 
 A slack matrix of the $8$-gon is given by  
\footnotesize
\begin{equation*} 
S_8=\begin{pmatrix}
0          & 1          & 1+\sqrt{2} & 2+\sqrt{2} & 2+\sqrt{2} & 1+\sqrt{2} & 1          & 0          \\
0          & 0          & 1          & 1+\sqrt{2} & 2+\sqrt{2} & 2+\sqrt{2} & 1+\sqrt{2} & 1          \\
1          & 0          & 0          & 1          & 1+\sqrt{2} & 2+\sqrt{2} & 2+\sqrt{2} & 1+\sqrt{2} \\
1+\sqrt{2} & 1          & 0          & 0          & 1          & 1+\sqrt{2} & 2+\sqrt{2} & 2+\sqrt{2} \\
2+\sqrt{2} & 1+\sqrt{2} & 1          & 0          & 0          & 1          & 1+\sqrt{2} & 2+\sqrt{2} \\
2+\sqrt{2} & 2+\sqrt{2} & 1+\sqrt{2} & 1          & 0          & 0          & 1          & 1+\sqrt{2} \\
1+\sqrt{2} & 2+\sqrt{2} & 2+\sqrt{2} & 1+\sqrt{2} & 1          & 0          & 0          & 1          \\
1          & 1+\sqrt{2} & 2+\sqrt{2} & 2+\sqrt{2} & 1+\sqrt{2} & 1          & 0          & 0          
\end{pmatrix}.
\end{equation*}
\normalsize

Let $\alpha_1=\sqrt{1+\sqrt{2}}$, $\alpha_2=\sqrt{2+\sqrt{2}}$, $\alpha_3=\frac{1}{\alpha_1}-\alpha_1$, $\alpha_4=\sqrt{\sqrt{2}}$ and $\alpha_5=\sqrt{1-\frac{1}{\alpha_1^2}}$.
A $\mathcal{S}^4_+$-factorization of $S_8$ is given by the following factors:

\small
\begin{equation*}
a^i=\left\{\begin{pmatrix}
1 \\
0 \\
0 \\
0
\end{pmatrix},
\begin{pmatrix}
0 \\
1 \\
0 \\
0
\end{pmatrix},
\begin{pmatrix}
0 \\
0 \\
1 \\
0
\end{pmatrix},
\begin{pmatrix}
1 \\
-\alpha_1 \\
-\alpha_1 \\
0
\end{pmatrix},
\begin{pmatrix}
1 \\
\alpha_3 \\
\alpha_3 \\
\frac{-1}{\alpha_1}
\end{pmatrix},
\begin{pmatrix}
0 \\
-1 \\
-2 \\
1
\end{pmatrix},
\begin{pmatrix}
0 \\
0 \\
1 \\
-1
\end{pmatrix},
\begin{pmatrix}
-1 \\
\frac{-1}{\alpha_1} \\
\frac{-1}{\alpha_1} \\
\frac{1}{\alpha_1}
\end{pmatrix}\right\},
\end{equation*}

\scriptsize
\begin{equation*}
b^j=\left\{\begin{pmatrix}
0 & 0 \\
0 & 0 \\
-1 & 0 \\
-1 & \alpha_1
\end{pmatrix},
\begin{pmatrix}
-1 & 0 \\
0 & 0 \\
0 & 0 \\
0 & -\alpha_2
\end{pmatrix},
\begin{pmatrix}
\alpha_1 & 0 \\
1 & 0 \\
0 & 0 \\
1 & \alpha_1
\end{pmatrix},
\begin{pmatrix}
-\alpha_2 & 0 \\
-\alpha_4 & 1 \\
0 & -1 \\
-\alpha_4 & 0
\end{pmatrix},
\begin{pmatrix}
\alpha_2 & 0 \\
0 & \alpha_2 \\
\frac{1}{\alpha_4} & \frac{-\alpha_1^2}{\alpha_2} \\
\frac{2}{\alpha_4} & \frac{-\sqrt{2}}{2}
\end{pmatrix},
\begin{pmatrix}
\alpha_1 & 0 \\
-\sqrt{2} & \alpha_4 \\
\sqrt{2} & -\alpha_4 \\
\sqrt{2} & -\alpha_4
\end{pmatrix},
\begin{pmatrix}
-1 & 0 \\
\alpha_1 & 0 \\
-\alpha_1 & 1 \\
-\alpha_1 & 1
\end{pmatrix},
\begin{pmatrix}
0 & 0 \\
\frac{-1}{\alpha_1} & \alpha_5 \\
\alpha_1 & 0 \\
-\alpha_3 & \alpha_5 
\end{pmatrix}\right\}.
\end{equation*}
\normalsize
\end{exmp}

\begin{exmp}
A slack matrix of the $10$-gon with $\phi=\frac{1+\sqrt{5}}{2}$ is given by   
\small
\begin{equation*} 
S_{10}=\begin{pmatrix}
0         & 0         & \phi^{-2} & 1         & \phi      & 2         & 2         & \phi      & 1         & \phi^{-2} \\
\phi^{-2} & 0         & 0         & \phi^{-2} & 1         & \phi      & 2         & 2         & \phi      & 1         \\
1         & \phi^{-2} & 0         & 0         & \phi^{-2} & 1         & \phi      & 2         & 2         & \phi      \\
\phi      & 1         & \phi^{-2} & 0         & 0         & \phi^{-2} & 1         & \phi      & 2         & 2         \\
2         & \phi      & 1         & \phi^{-2} & 0         & 0         & \phi^{-2} & 1         & \phi      & 2         \\
2         & 2         & \phi      & 1         & \phi^{-2} & 0         & 0         & \phi^{-2} & 1         & \phi      \\
\phi      & 2         & 2         & \phi      & 1         & \phi^{-2} & 0         & 0         & \phi^{-2} & 1         \\
1         & \phi      & 2         & 2         & \phi      & 1         & \phi^{-2} & 0         & 0         & \phi^{-2} \\
\phi^{-2} & 1         & \phi      & 2         & 2         & \phi      & 1         & \phi^{-2} & 0         & 0         \\
0         & \phi^{-2} & 1         & \phi      & 2         & 2         & \phi      & 1         & \phi^{-2} & 0         \\
\end{pmatrix}.
\end{equation*}

Let $\alpha_1=(\sqrt{2}\phi)^{-1/2}$, $\alpha_2=(\sqrt{2}/\phi)^{1/2}$, $\alpha_3=\sqrt{2/\phi}$, $\alpha_4=\sqrt{\sqrt{2}}\phi$, $\alpha_5=-\phi^{3/2}$ and $\alpha_6=\sqrt{\sqrt{5}-1}$.
A $\mathcal{S}^5_+$-factorization of $S_{10}$ is given by the following factors:
\tiny 
\begin{align*}
  a^i=\left\{\begin{pmatrix}
0 & \alpha_1^{-1} & 0\\
0 & \alpha_4 & 0 \\
1 & \alpha_5 & 0 \\
0 & -1 & 0 \\
0 & 0 & 1
\end{pmatrix},
\begin{pmatrix}
0 & (\alpha_1\phi)^{-1} & 0\\
\alpha_2 & \alpha_4\phi^{-1} & 0 \\
-1 & \alpha_5\phi^{-1} & 0 \\
0 & -\phi^{-1} & 0 \\
0 & 0 & \sqrt{\phi}
\end{pmatrix},
\begin{pmatrix}
0 & 0 & 0\\
\alpha_2 & 0 & 0 \\
-1 & 0 & 0 \\
0 & \phi^{-1} & 0 \\
0 & 0 & \sqrt{\phi}
\end{pmatrix},
\begin{pmatrix}
0 & 0 & 0\\
0 & 0 & 0 \\
1 & 0 & 0 \\
0 & 1 & 0 \\
0 & 0 & 1
\end{pmatrix},
\begin{pmatrix}
\alpha_2 & 0 & 0 \\
0 & 0 & 0  \\
0 & \phi & 0  \\
0 & 0 & \alpha_3  \\
0 & 0 & 0 
\end{pmatrix},
\begin{pmatrix}
0 & 0 & \alpha_1^{-1} \\
0 & 0 & \alpha_2  \\
1 & 0 & \sqrt{2\phi}-1  \\
0 & 1 & 0  \\
0 & 0 & -1 
\end{pmatrix}\right.,
  \mspace{50mu}
  \notag\\
     \left.\begin{pmatrix}
0 & 0 & \alpha_4 \\
\alpha_2 & 0 & 2^{1/4}  \\
-1 & 0 & \sqrt{2}\phi-\sqrt{\phi}  \\
0 & \phi^{-1} & 0  \\
0 & 0 & -\sqrt{\phi} 
\end{pmatrix},
\begin{pmatrix}
0 & (\alpha_1\phi)^{-1} & \alpha_4\\
\alpha_2 & \alpha_4\phi^{-1} & 2^{1/4} \\
-1 & \alpha_5\phi^{-1} & \sqrt{2}\phi-\sqrt{\phi} \\
0 & -\phi^{-1} & 0 \\
0 & 0 & -\sqrt{\phi}
\end{pmatrix},
\begin{pmatrix}
0 & \alpha_1^{-1} & \alpha_1^{-1}\\
0 & \alpha_4 & \alpha_2 \\
1 & \alpha_5 & \sqrt{2\phi}-1 \\
0 & -1 & 0 \\
0 & 0 & -1
\end{pmatrix},
\begin{pmatrix}
\alpha_2 & 0 & \alpha_3\phi^{-1}\\
0 & 0 & \alpha_4\alpha_3 \\
0 & \phi & \alpha_5\alpha_3 \\
0 & 0 & -\alpha_3 \\
0 & 0 & 0
\end{pmatrix}\right\},
\end{align*}

\footnotesize
\begin{equation*}
b^j=\left\{\begin{pmatrix}
0 \\
\alpha_1 \\
0 \\
\sqrt{\phi} \\
0
\end{pmatrix},
\begin{pmatrix}
\alpha_1 \\
0 \\
0 \\
1 \\
0
\end{pmatrix},
\begin{pmatrix}
0 \\
\alpha_1 \\
\phi^{-1} \\
0 \\
0
\end{pmatrix},
\begin{pmatrix}
\alpha_1 \\
0 \\
0 \\
0 \\
0
\end{pmatrix},
\begin{pmatrix}
0 \\
\alpha_1 \\
0 \\
0 \\
0
\end{pmatrix},
\begin{pmatrix}
0 \\
\alpha_1 \\
0 \\
0 \\
\phi^{-1}
\end{pmatrix},
\begin{pmatrix}
\alpha_1 \\
0 \\
0 \\
0 \\
1
\end{pmatrix},
\begin{pmatrix}
0 \\
\alpha_1 \\
\phi^{-1} \\
0 \\
\alpha_6
\end{pmatrix},
\begin{pmatrix}
\alpha_1 \\
0 \\
0 \\
1 \\
1
\end{pmatrix},
\begin{pmatrix}
0 \\
\alpha_1 \\
0 \\
\sqrt{\phi} \\
\phi^{-1}
\end{pmatrix}\right\}.
\end{equation*}

\normalsize
\end{exmp}

\subsection{Adaptation for Completely PSD matrices} \label{completelypsd}

For NMF (\ref{nmf-opt}) involving a symmetric $n$-by-$n$ matrix $X$, 
the additional constraint requiring $W$ and $H$ to be equal to each other  
leads to an optimization problem known as symmetric NMF (SymNMF).  
Specific numerical algorithms have been designed for this problem having applications in data mining \cite{ho2008nonnegative,kuang2015symnmf,vandaele2016efficient}.  
When an exact factorization is possible, that is, $X=HH^T$ for a nonnegative $n$-by-$k$ matrix $H$, the matrix $X$ is said to be completely positive. 
The smallest integer $k$ for which such an exact factorization exists is referred as the cp-rank of $X$~\cite{berman2003completely}. 

By analogy with completely positive matrices, a completely positive semidefinite matrix $X$ is defined as a $n$-by-$n$ symmetric matrix for which there exists a set $A^1,...,A^n \in \mathcal{S}^k_+$ such that $X_{i,j}=\left\langle A^i,A^j\right\rangle$.
The smallest integer $k$ for which it is possible to write such a factorization is called the cpsd-rank of $X$; see, e.g.,~\cite{gribling2017matrices,prakash2016completely}. 
As opposed to problem (\ref{psd-obj-sub}), the symmetric version 
\begin{equation}\label{cpsd-obj-sub}
\min_{\substack{A^i\in\mathcal{S}_+^{k} \\ i=1,...,n}}\sum_{i=1}^n\sum_{j=1}^n\left(X_{i,j}-\left\langle A^i,A^j\right\rangle\right)^2,
\end{equation}
is no longer convex even when all factors $A^i$'s are fixed but one. 
However, it is possible to adapt the methods developed in Section \ref{sec-cd} in order to handle (\ref{cpsd-obj-sub}).
We propose to keep the problem with two sets of variables but we add a penalty term to (\ref{psd-obj-newvariables}) with a scalar $\gamma>0$ in order to enforce the similarity between $a^i$ and $b^i$ for $i=1,...,n$, similarly as done for Symmetric NMF in~\cite{ho2008nonnegative,kuang2015symnmf}:  
\begin{equation*} 
\min_{\substack{a^i\in\mathbb{R}^{k\times r_i} \\ i=1,...,n}}f=\sum_{i=1}^n\sum_{j=1}^n\left(X_{i,j}-\sum_{h=1}^{r_i}\sum_{l=1}^ {r_j}\left({a_{:,h}^i}^Tb_{:,l}^j\right)^2\right)^2 + \gamma\sum_{i=1}^n\|a^i-b^i\|_2^2.
\end{equation*} 
This modification of the objective function has limited consequences on Algorithms~\ref{algcyclicCD} and~\ref{algGSCD} since the additional terms are quadratic. The entry of the gradient corresponding to the variable $a_{p,q}^i$ is given by 
\begin{equation*} 
\nabla_{a_{p,q}^i}f = c_3{a_{p,q}^i}^3+c_2{a_{p,q}^i}^2+(c_1+2\gamma){a_{p,q}^i}+(c_0-2\gamma b^i_{p,q}).
\end{equation*}
With this change, we are now able to compute symmetric factorizations.
\begin{exmp} 
	The symmetric $6$-by-$6$ matrix  
	\begin{equation*} 
	P_4=\begin{pmatrix}
	2 & 1 & 1 & 1 & 1 & 0\\
	1 & 2 & 1 & 1 & 0 & 1\\
	1 & 1 & 2 & 0 & 1 & 1\\
	1 & 1 & 0 & 2 & 1 & 1\\
	1 & 0 & 1 & 1 & 2 & 1\\
	0 & 1 & 1 & 1 & 1 & 2
	\end{pmatrix}, 
	\end{equation*}
	 as defined in Section \ref{psd-datasets} has a symmetric factorization with $k=4$ 	with the factors 
	\begin{equation*}
	a^i=\left\{\begin{pmatrix}
	1 & 0 \\
	0 & 1 \\
	0 & 0 \\
	0 & 0
	\end{pmatrix},
	\begin{pmatrix}
	1 & 0 \\
	0 & 0 \\
	0 & 1 \\
	0 & 0
	\end{pmatrix},
	\begin{pmatrix}
	1 & 0 \\
	0 & 0 \\
	0 & 0 \\
	0 & 1
	\end{pmatrix},
	\begin{pmatrix}
	0 & 0 \\
	0 & 1 \\
	1 & 0 \\
	0 & 0
	\end{pmatrix},
	\begin{pmatrix}
	0 & 0 \\
	1 & 0 \\
	0 & 0 \\
	0 & 1
	\end{pmatrix},
	\begin{pmatrix}
	0 & 0 \\
	0 & 0 \\
	1 & 0 \\
	0 & 1
	\end{pmatrix}\right\}.
	\end{equation*}
\end{exmp}

The choice of the parameter $\gamma$ can be made in different ways and should be increased in the course of the optimization process in order to ensure that $a^i$ converges to $b^i$ for all $i$.
For this particular example, we simply used $\gamma=1$ which gave us the desired result.

\subsection{Adaptation for the square root rank} \label{square-root-rank}

Given a nonnegative matrix $X$, a \textit{Hadamard square root} of $X$ is defined as a matrix obtained by replacing the $(i,j)$th entry of $X$ by either $\sqrt{X_{i,j}}$, or $-\sqrt{X_{i,j}}$.
Hence there are $2^{N}$ possible Hadamard square root matrices for a matrix $X$ with $N$ non-zero entries. 
The \textit{square root rank} of a nonnegative matrix $X$ is defined as the minimum rank among the ranks of all the Hadamard square root of $X$. If a nonnegative matrix $X$ has square root rank $k$, then there is an exact PSD factorization of $X$ with rank-1 factors of size $k$; 
see Proposition 6.2. in \cite{FGP14} (hence the square root rank of $X$ is an upper bound on the psd-rank of $X$).   
Therefore, we can use Algorithms~\ref{algcyclicCD} and~\ref{algGSCD} with $r_i=1$ for all $i$ to try to compute the square root rank of $X$. Note that computing this quantity is NP-hard as well~\cite{FGP14}.  

\begin{exmp}
For the 8-gon (and its slack matrix $S_8$; see Example~\ref{ex8gon}), we have computed such a rank-one decomposition with $k=6$ and $r_i = 1$ for all $i$. 
Note that to compute this decomposition, we had to use many different starting points (around a thousand) and manually fix some entries 
of the $a^i$'s and $b^j$'s to zero.  
In order to present this decomposition, let us define
\small
\[
S=\begin{pmatrix}
     0  &  -1  & -1  &  1  &  -1  &  -1  &  1  &  0 \\
     0  &   0  &  1  &  1  &   1  &  -1  &  1  & -1 \\
     1  &   0  &  0  & -1  &  -1  &   1  & -1  &  1 \\
     1  &  -1  &  0  &  0  &  -1  &  -1  &  1  & -1 \\
     1  &  -1  &  1  &  0  &   0  &   1  &  1  &  1 \\
     1  &   1  &  1  & -1  &   0  &   0  & -1  & -1 \\
     1  &   1  &  1  &  1  &   1  &   0  &  0  &  1 \\
    -1  &   1  & -1  &  1  &  -1  &  -1  &  0  &  0 \\
\end{pmatrix},
\]
\[
W= \begin{pmatrix}
     0         &  -1        & -\alpha_1  &  \alpha_2  &  -\alpha_2  &  -\alpha_1 \\
     0         &   0        &  1         &  \alpha_1  &   \alpha_2  &  -\alpha_2 \\
     1         &   0        &  0         &  -1        &  -\alpha_1  &  \alpha_2  \\
     \alpha_1  &  -1        &  0         &   0        &  -1         &  -\alpha_1 \\
     \alpha_2  & -\alpha_1  &  1         &   0        &   0         &   1        \\
     \alpha_2  &  \alpha_2  &  \alpha_1  &  -1        &   0         &   0        \\
     \alpha_1  &  \alpha_2  &  \alpha_2  &  \alpha_1  &   1         &   0        \\
     -1        &  \alpha_1  & -\alpha_2  &  \alpha_2  &  -\alpha_1  &  -1        \\
\end{pmatrix}, \text{ and }
H= \begin{pmatrix}
     1  &  0  &  0  &  0  &  0  &  0  &  \frac{\alpha_2}{\alpha_1-1}    &  0       \\
     0  &  1  &  0  &  0  &  0  &  0  &  0                              &  \frac{1-\alpha_1}{\alpha_2} \\
     0  &  0  &  1  &  0  &  0  &  0  &  \frac{1+\alpha_1}{1-\alpha_1}  &  0       \\
     0  &  0  &  0  &  1  &  0  &  0  &  0                              &  1       \\
     0  &  0  &  0  &  0  &  1  &  0  &  \frac{\alpha_2}{\alpha_1-1}    &  0       \\
     0  &  0  &  0  &  0  &  0  &  1  &  0                              &  \frac{1+\alpha_1}{\alpha_2} \\
\end{pmatrix},
\]
with $\alpha_1=\sqrt{1+\sqrt{2}}$ and $\alpha_2=\sqrt{2+\sqrt{2}}$. 
Denoting $\sqrt[+]{S_8}$ the nonnegative Hadamard square root of $S_8$, one can check that $S\circ\sqrt[+]{S_8}=WH$ implying that the square root rank of $S_8$ is at most 6. 
Our algorithms were not able to compute such a decomposition for $k=5$ (relative error always at least $0.6\%$). 

\normalsize
\end{exmp}


\section{Conclusion} \label{sec-conclusion} 

In this work, we introduced different algorithms for solving numerically the PSD factorization problem (\ref{psd-opt}).
These algorithms are based on an alternating strategy in order to solve convex subproblems.  
The first method proposed uses PSD matrices as variables and implements a fast projected gradient method. 
The second idea is to apply the coordinate descent (CD) framework after having expressed the problem as an unconstrained optimization problem.
Numerical experiments have been conducted to assess the performances of the different methods, and we observed that CD with the Gauss-Southwell rule performs consistently the best. 
Finally, we have illustrated the ability of our algorithms to help in the computation of non-trivial factorizations for regular $n$-gons, for symmetric PSD factorizations and for the square root rank.  
Note that an earlier version of our code was also used successfully in~\cite{kubjas2015positive}.  

An important direction for future research is the development of a globalization framework, such as in \cite{VGGT14} for NMF, 
in order to escape local minima and generate, in average, better solutions than with a simple multi-start strategy as used in this paper.

	\small 

\bibliographystyle{plain}
\bibliography{references}

\normalsize 

\appendix
\section{Minimizer of a quartic polynomial}  \label{appA}

This appendix is devoted to the description of the algorithm for computing the minimum value of the univariate quartic polynomials arising in Section \ref{sec-cd}, which are of the form
\[
	f(x)=c_3\frac{x^4}{4}+c_2\frac{x^3}{3}+c_1\frac{x^2}{2}+c_0x+K.
\]
Finding the minimizer of such a function can be done by enumerating the roots of $f'(x)$ and choosing the one minimizing $f$.
They are at most three roots to a cubic equation but we will show that in our case, we only have to consider two of them.
We use the well-known trick known under the name \textit{Cardano's method} (see for example \cite{cardano1968ars}) to reduce the identification of the solutions of $c_3x^3+c_2x^2+c_1x+c_0=0$ to the computation of the roots of the following \textit{depressed} cubic polynomial
\[
	p(t)=t^3+at+b,
\]
where $a=\frac{3c_3c_1-c_2^2}{3c_3^2}$, $b=\frac{2c_2^3-9c_3c_2c_1+27c_3^2c_0}{27c_3^3}$ and  the substitution $t-\frac{c_2}{3c_3}$ for $x$.
\begin{itemize}
\item When the quantity $\Delta=4a^3+27b^2$ is positive, there is only one real root which has a closed form expression:
\[
	t^{\ast}=\sqrt[3]{\frac{1}{2}\left(-b+\sqrt{\frac{\Delta}{27}}\right)}+\sqrt[3]{\frac{1}{2}\left(-b-\sqrt{\frac{\Delta}{27}}\right)}.
\]
\item When $\Delta$ is negative (which implies $a<0$), there are at most three real roots
\[
	t^{\ast}_l=r\cos\left(\theta + \frac{2l\pi}{3}\right)\text{ for }l=0,1,2,
\]
with $r=2\sqrt{\frac{-a}{3}}$ and $\theta=\frac{1}{3}\arccos\left(\frac{3b}{2a}\sqrt{\frac{-3}{a}}\right)$.
It is easy to check that $t_1\leq t_2\leq t_0$ and since the coefficient of the leading term $c_3\geq 0$ in our case, the root $t_2$ always corresponds to a local maximum of $f(x)$.
\end{itemize}
Taking into accounts the previous observations, Algorithm~\ref{cardanomethod} is the pseudo-code of the method used in our CD schemes to determine in $O(1)$ operations the minimizer of the univariate quartic polynomials. 

\begin{algorithm} 
\caption{$x=CardanoMethod(c_3,c_2,c_1,c_0)$}
\label{cardanomethod} 
\begin{algorithmic}[1]
\STATE INPUT: $c_3 \in \mathbb{R}_0^+, c_2 \in \mathbb{R}, c_1 \in \mathbb{R}, c_0 \in \mathbb{R}$
\STATE OUTPUT: $x^{\ast}=\arg \min_{x}$ $c_3\frac{x^4}{4}+c_2\frac{x^3}{3}+c_1\frac{x^2}{2}+c_0x$.
\STATE $a = \frac{3c_3c_1-c_2^2}{3c_3^2}$
\STATE $b = \frac{2c_2^3-9c_3c_2c_1+27c_3^2c_0}{27c_3^3}$
\STATE $\Delta = 4a^3+27b^2$
\IF{$\Delta\leq 0$}
\STATE $t_0=2\sqrt{\frac{-a}{3}}\cos\left(\frac{1}{3}\arccos\left(\frac{3b}{2a}\sqrt{\frac{-3}{a}}\right)\right)$
\STATE $t_1=2\sqrt{\frac{-a}{3}}\cos\left(\frac{1}{3}\arccos\left(\frac{3b}{2a}\sqrt{\frac{-3}{a}}\right) + \frac{2\pi}{3}\right)$
\IF{$\frac{t_0^4}{4}+a\frac{t_0^2}{2}+bt_0 < \frac{t_1^4}{4}+a\frac{t_1^2}{2}+bt_1$}
\STATE $t^{\ast}=t_0$
\ELSE
\STATE $t^{\ast}=t_1$
\ENDIF
\ELSE
\STATE $t^{\ast} = \sqrt[3]{\frac{1}{2}\left(-b+\sqrt{\frac{\Delta}{27}}\right)}+\sqrt[3]{\frac{1}{2}\left(-b-\sqrt{\frac{\Delta}{27}}\right)}$
\ENDIF
\STATE $x^{\ast}=t^{\ast}-\frac{c_2}{3c_3}$
\end{algorithmic}
\end{algorithm}

\end{document}